\documentclass[11pt]{article}
\usepackage[top=1in, bottom=1in, left=1in, right=1in]{geometry}
\usepackage{tikz}
\usetikzlibrary{shapes,arrows,matrix,patterns,positioning}
\usepackage{color}
\usepackage{graphicx}
\usepackage{algorithmic}
\usepackage{amssymb}
\usepackage{epstopdf}
\usepackage[ruled]{algorithm2e}
\usepackage{float}
\usepackage{amsmath,amsthm,bm,color,epsfig,enumerate,caption}
\usepackage{hyperref}
\usepackage{booktabs}
\usepackage{multirow}
\usepackage{array}
\usepackage{bm}

\newtheorem{theorem}{Theorem}[section]

\newtheorem{assumption}[theorem]{Assumption}

\newtheorem{algo}[theorem]{Algorithm}

\renewcommand{\appendix}[1]{
\section*{Appendix: #1}
}

\renewcommand{\O}{\mathcal{O}}
\renewcommand{\O}{O}

\newcommand{\bbC}{\mathbb{C}}
\newcommand{\bbR}{\mathbb{R}}
\newcommand{\R}{\mathbb{R}}
\newcommand{\eps}{\epsilon}

\newcommand{\diag}{\text{ diag}}
\newcommand{\coef}[2]{\lambda^{#1}_{#2}}
\newcommand{\Coef}[2]{\Lambda^{#1}_{#2}}
\newcommand{\tCoef}[2]{\tilde{\Lambda}^{#1}_{#2}}

\makeatletter
\newcommand*{\extendadd}{
  \mathbin{
    \mathpalette\extend@add{}
  }
}
\newcommand*{\extend@add}[2]{
  \ooalign{
    $\m@th#1\leftrightarrow$%
    \vphantom{$\m@th#1\updownarrow$}
    \cr
    \hfil$\m@th#1\updownarrow$\hfil
  }
}
\makeatother

\begin{document}

\title{Interpolative Butterfly Factorization}

\author{Yingzhou Li$^\sharp$,
    Haizhao Yang$^\dagger$,
  \vspace{0.1in}\\
  $\sharp$ ICME, Stanford University\\
  $\dagger$ Department of Mathematics, Duke University
}

\maketitle

\begin{abstract}
This paper introduces the interpolative butterfly
factorization for nearly optimal implementation of several transforms in
harmonic analysis, when their explicit formulas satisfy certain analytic
properties and the matrix representations of these transforms
satisfy a complementary low-rank property.
A preliminary interpolative butterfly factorization is constructed based on
interpolative low-rank approximations of the complementary low-rank matrix.
A novel sweeping matrix compression technique further compresses
the preliminary interpolative butterfly factorization
via a sequence of structure-preserving low-rank
approximations. The sweeping procedure propagates the low-rank property among
neighboring matrix factors to compress dense submatrices
in the preliminary butterfly factorization to obtain an optimal one
in the butterfly scheme.
For an $N\times N$ matrix, it takes $O(N\log N)$ operations
and complexity to construct the factorization as a product of
$O(\log N)$ sparse matrices, each with $O(N)$ nonzero entries. Hence,
it can be applied rapidly in $O(N\log N)$ operations. Numerical results
are provided to demonstrate the effectiveness of this algorithm.
\end{abstract}



\section{Introduction} \label{sec:intro} One key problem in computational
harmonic analysis is the rapid evaluation of various transforms to make
large-scale scientific computation feasible. These transforms are
essentially matrix-vector multiplications, $u=Kg$, where the kernel matrix
$K$  is the discrete representation of a transform and the vector $g$ is the
discrete representation of a function to be transformed.  Inspired by the
idea of the butterfly algorithm initially proposed in \cite{Butterfly1} and
later extended in \cite{Butterfly2}, the recently proposed butterfly
factorization \cite{BF,MBF} factorizes a complementary low-rank matrix $K$
of size $N\times N$ into a product of $O(\log N)$ sparse matrices, each with
$O(N)$ nonzero entries. After factorization, the application of
$K$ has nearly optimal\footnote{Through out the paper, the ``nearly optimal"
refers to the nearly optimal constant in the complexity.} operation and
memory complexity of order $N\log N$.  Since a wide range of transforms in
harmonic analysis admits a matrix representation satisfying the
complementary low-rank property
\cite{FIO09,FIO12,FIO14,SHT,Butterfly2,FIO13}, the butterfly factorization,
once constructed, is a nearly optimal fast algorithm to evaluate these
transforms.  However, the construction of the butterfly factorization
requires $O(N^2)$ operations in \cite{Butterfly2} and requires $O(N^{1.5})$
operations in \cite{BF,MBF}, which might still be too expensive in real
applications.  This paper introduces the interpolative butterfly
factorization to construct the factorization in $O(N\log N)$ operation and
memory complexity, if the continuous kernel ${K}$ is explicitly available. 
The interpolative butterfly factorization is a combination of the butterfly algorithm in 
\cite{FIO09} and a novel structure-preserving matrix compression technique. 
Hence, the proposed method can be considered as an 
optimized sparse matrix representation of the butterfly algorithm in \cite{FIO09} with 
a smaller prefactor in the operation complexity.

Another key problem in modern large-scale computation is
the parallel scalability
of fast algorithms. Although there have been various fast algorithms like
 the (nonuniform) fast Fourier transform (FFT)
\cite{PolarFFT,FFT,NUFFT}, the FFT in polar and spherical coordinates
\cite{Knockaert,SHT2,Alex,SHT},
the parallel scalability of some of
these traditional algorithms might be still limited
in high performance computing. This motivates much effort to
improve their parallel scalability \cite{PFFT2,PFFT1,wavemoth,ying20053d}.
For the same purpose, the interpolative butterfly factorization
is proposed as a general framework for highly parallel scalable
implementation of a wide range of transforms in harmonic analysis.
Since the construction of the butterfly factorization
is just a few essentially independent low-rank approximations and
the application is a sequence of small matrix-vector multiplications,
the butterfly factorization framework significantly reduces communication
if implemented in parallel computation.

To be more specific, the interpolative butterfly factorization is proposed
for the rapid application of integral transforms of the form
\begin{equation}
  \label{eqn:FIO}
  u(x) = \int_{\mathbb{R}^d}a(x,\xi)  e^{2\pi \i \Phi(x,\xi)}g(\xi) d\xi,
\end{equation}
where $d$ is the dimension,
and ${K}(x,\xi)=a(x,\xi)  e^{2\pi \i \Phi(x,\xi)}$ is the kernel
function that satisfies following properties:
\begin{assumption}{Smoothness properties}
\label{asm}

\begin{itemize}
\item $a(x,\xi)$ is an amplitude function that is smooth both in
  $x$ and $\xi$;
\item $\Phi(x,\xi)$ is a phase function that is real analytic for
  $x$ and $\xi$ and obeys the homogeneity condition of degree $1$ in
  $\xi$, namely, $\Phi(x,\lambda \xi)=\lambda\Phi(x,\xi)$ for
  $\lambda>0$.
\end{itemize}
\end{assumption}
This transform is also known as the Fourier integral operator (FIO).
FIOs are a wide class of operators in harmonic analysis including
 the (nonuniform) Fourier transform, pseudo-differential operators,
 and the generalized Radon transform. All these are popular tools in
 computational physics and chemistry \cite{Yingwave,Hu,PFFT,invRadon,invFIO},
 imaging science \cite{PR1,Imaging1,SP2,Zhizhen:2014}.
 For higher dimensional FIOs, the phase
function might not be smooth when $\xi= 0$. Fortunately, the
interpolative butterfly factorization can be adapted to this case
following the idea in the multiscale butterfly algorithm in \cite{FIO14}.

In most examples, since $a(x,\xi)$ is a smooth symbol of order zero
and type $(1,0)$ \cite{symbol1,FIO07,symbol2,symbol3},
$a(x,\xi)$ is numerically low-rank in the joint $X$ and $\Omega$
domain and its numerical treatment is relatively easy. Therefore, we
will simplify the problem by assuming $a(x,\xi)=1$ in the following
 discussion.

 In a typical setting, it is often assumed that the function $g(\xi)$
 decays sufficiently fast so that one can embed the problem in a
 sufficiently large periodic cell. Without loss of generality, a simple
 discretization considers functions given on a
 Cartesian grid
\begin{equation}
  \label{eqn:X}
  X = \left\{ x =\left( \frac{n_1}{N^{1/d}},\dots, \frac{n_d}{N^{1/d}}\right), 0 \leq  n_1,\dots,n_d < N^{1/d}\text{ with }
  n_1,\dots,n_d \in \mathbb{Z} \right\}
\end{equation}
in a unit box in $x$ and defines the discrete integral transform by
\begin{equation}
\label{eqn:DFIO}
u(x) = \sum_{\xi\in\Omega} {K}(x,\xi) g(\xi), \quad x\in X,
\end{equation}
where
\begin{equation}
  \label{eqn:Omega}
  \Omega = \left\{ \xi =(n_1,\dots,n_d),- \frac{N^{1/d}}{2} \leq n_1,\dots,n_d < \frac{N^{1/d}}{2}\text{ with }
  n_1,\dots,n_d \in \mathbb{Z} \right\}.
\end{equation}
Using the notation in numerical linear algebra, the evaluation
of \eqref{eqn:DFIO} is a matrix-vector multiplication $u=Kg$.
Under Assumption \ref{asm}, it
can be proved that $K$ is essentially complementary low-rank \cite{FIO09,FIO14}.

\subsection{Complementary low-rank matrices and interpolative butterfly factorization}

\begin{figure}[ht!]
  \begin{center}
    \begin{tabular}{ccccc}
      \includegraphics[height=1.1in]{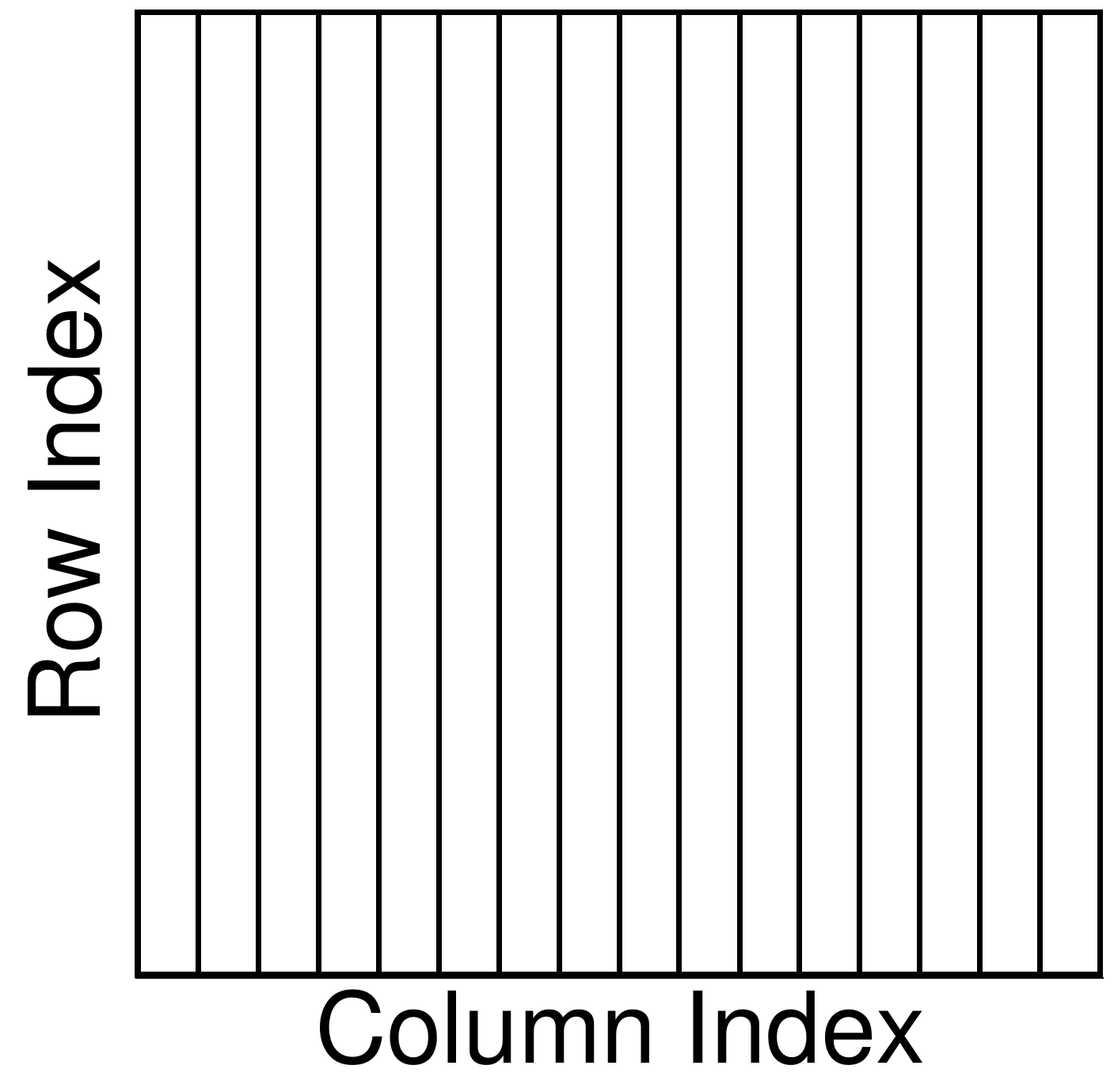}&
      \includegraphics[height=1.1in]{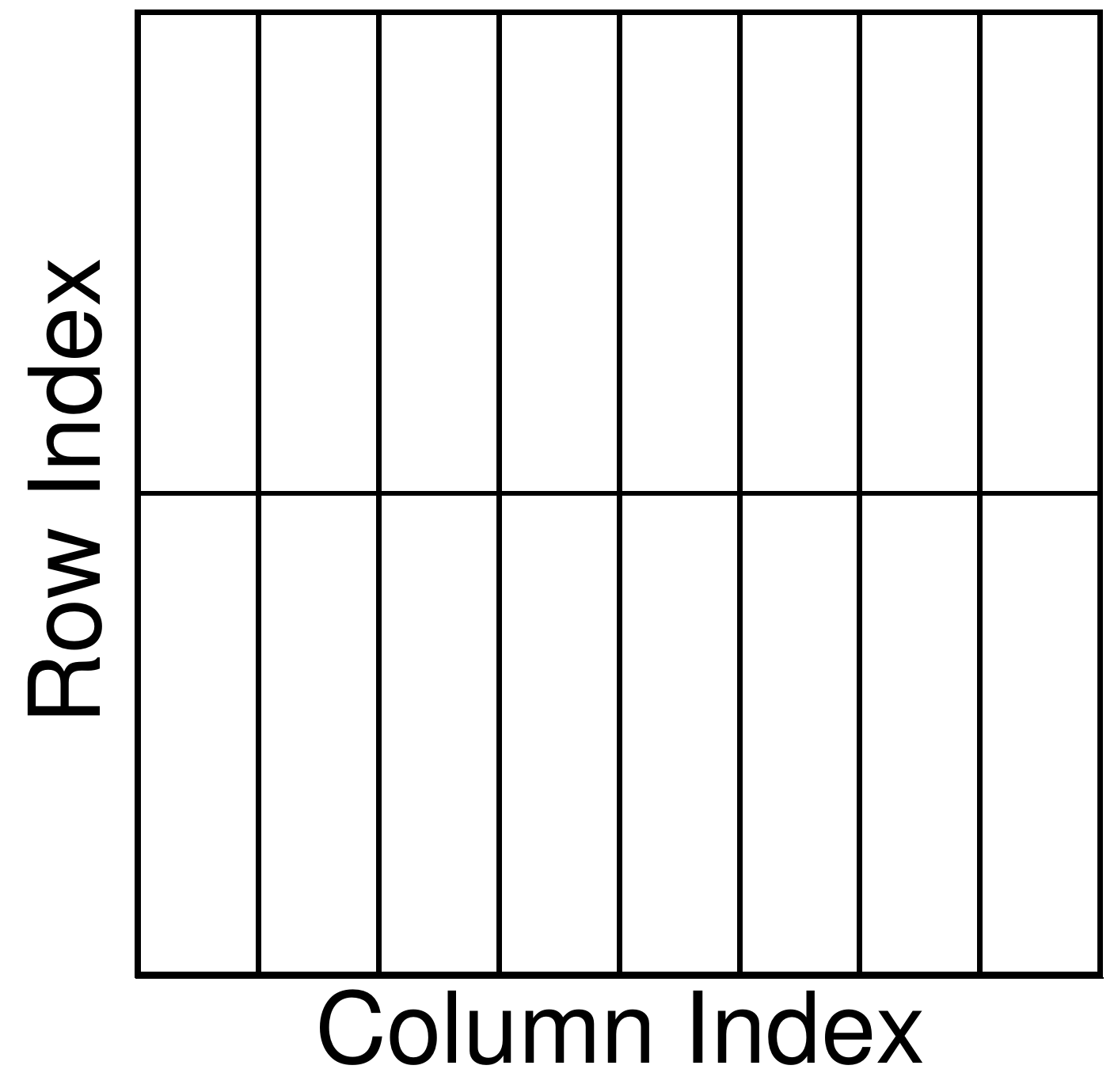}&
      \includegraphics[height=1.1in]{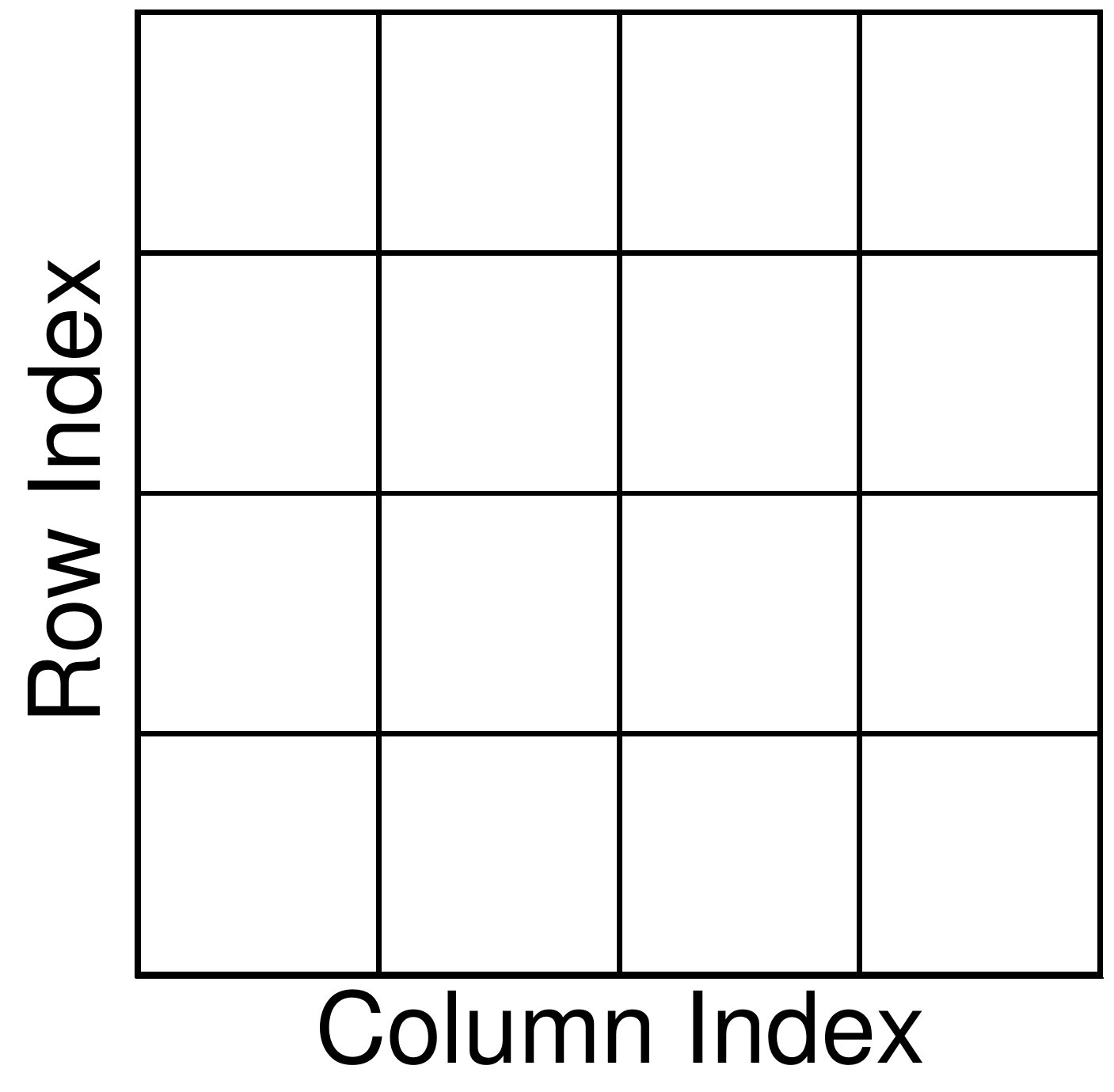}&
      \includegraphics[height=1.1in]{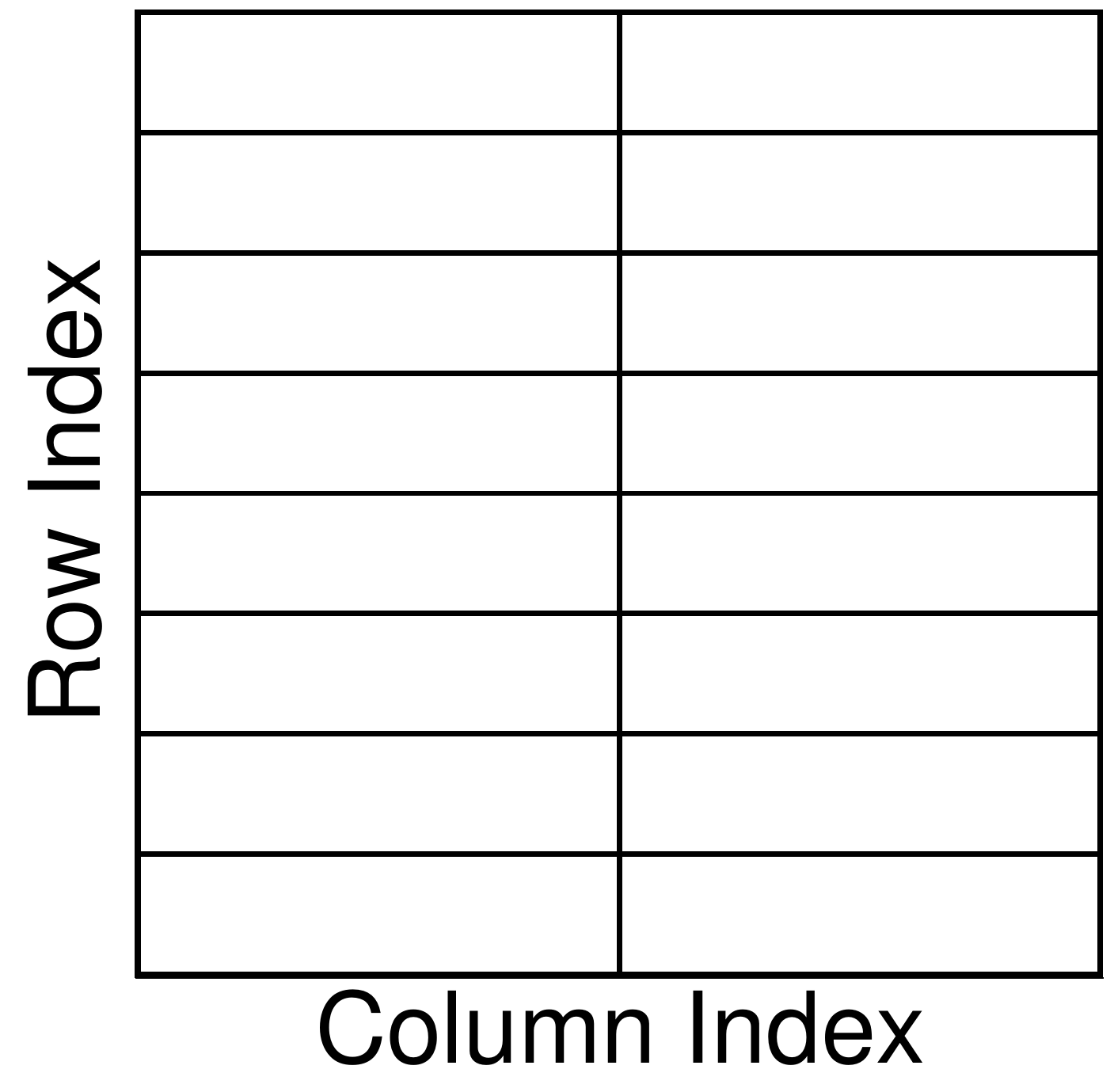}&
      \includegraphics[height=1.1in]{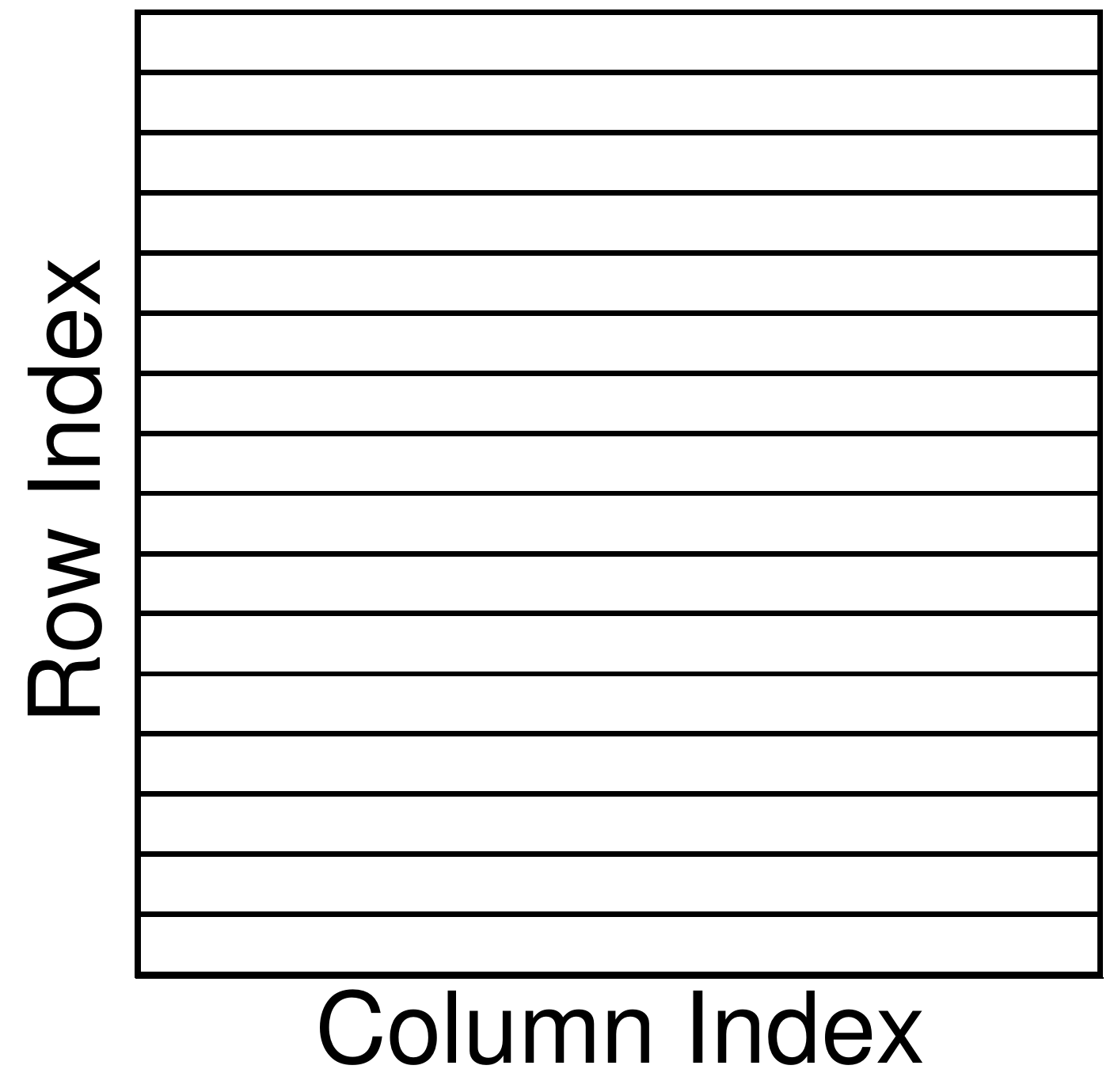}
    \end{tabular}
  \end{center}
  \caption{Hierarchical decomposition of the row and column indices
    of a one-dimensional complementary low-rank matrix of size
     $16\times 16$.  The trees $T_{X}$ ($T_{\Omega}$)
    has a root containing $16$ column (row) indices and leaves
    containing a single column (row) index.  The rectangles above
    indicate some of the low-rank submatrices.}
\label{fig:submatrices}
\end{figure}

Complementary low-rank matrices have been widely studied in \cite{Butterfly3,Butterfly4,BF,MBF,Butterfly2,Butterfly1,
Butterfly5}.
 Let $X$ and $\Omega$ be point sets (not necessary uniformly distributed) in $\bbR^d$ for
some dimension $d$. When the kernel $K(x,\xi)$ is discretized on $X\times \Omega$,
points in $X$ and $\Omega$ are indexed with row and column indices
in the matrix $K$. For simplicity, we also use
$X$ and $\Omega$ to denote the sets of row and column indices.
Two trees $T_X$ and $T_\Omega$ of the same depth $L=O(\log N)$,
associated with $X$ and $\Omega$ respectively,
are constructed by dyadic partitioning. Denote the root level of the tree as level $0$ and
the leaf one as level $L$. Such a matrix $K$ of size $O(N)\times O(N)$
is said to satisfy the {\bf complementary low-rank property} if for
any level $\ell$, any node $A$ in $T_X$ at level $\ell$, and any node
$B$ in $T_\Omega$ at level $L-\ell$, the submatrix $K_{A,B}$, obtained
by restricting $K$ to the rows indexed by the points in $A$ and the
columns indexed by the points in $B$, is numerically low-rank, i.e.,
for a given precision $\epsilon$ there exists a low-rank approximation
of $K_{A,B}$ with an error bounded by $\epsilon$ and the
rank bounded polynomially in $\log(1/\epsilon)$ and is independent of $N$.
See Figure \ref{fig:submatrices} for an illustration of low-rank submatrices
in a complementary low-rank matrix of size $16\times 16$.

It will be shown that, for a complementary low-rank matrix $K$,
the matrix-vector multiplication $u=Kg$
can be carried out efficiently via a preliminary
{\bf interpolative butterfly factorization} (IBF) constructed by
interpolative low-rank approximations:
\[
K \approx  U^L G^{L-1}\cdots G^{h}  M^{h} (H^{h})^* \dots (H^{1})^* (V^0)^*,
\]
where the depth $L=\O(\log N)$ of $T_X$ and $T_\Omega$ is assumed to be even,
$h=L/2$ is a middle level index, all factors are sparse matrices with
$\O(N)$ nonzero entries and a large prefactor.
Dense blocks in these sparse factors
come from interpolative low-rank approximations of low-rank submatrices
as illustrated in Figure \ref{fig:submatrices}.
When the kernel function ${K}$ satisfies Assumption \ref{asm}, each low-rank
approximation of these submatrices can be constructed explicitly by
Lagrange interpolation with $q$ Chebyshev grid points in each dimension,
resulting in dense submatrices of size $q^d\times q^d$ in sparse factors.
Hence, all the matrix factors
are available explicitly. However, the low-rank
approximation by interpolation is not optimal in the sense that
it over estimates the numerical rank $r_0$ of the low-rank matrix, i.e. $q^d> r_0$. Hence, dense
blocks in these sparse factors
can be further compressed to smaller submatrices of size $r_0\times r_0$
 by a truncated SVD.

The preliminary IBF above can be further compressed by a novel sweeping
matrix compression method in two stages. In the sweep-out stage, a
sequence of structure-preserving matrix compression is conducted:
\[
{M^h}\approx {C^h}{\bar{M}^h}{(R^h)}^*,
\]
\begin{equation*}
  {G^\ell}  {C^\ell} \approx
 {C^{\ell+1}} {\bar{G}^{\ell}}
\end{equation*}
for $\ell=h,h+1,\dots,L-1$,
\begin{equation*}
( {H^\ell}  {R^{\ell}})^* \approx ( {\bar{H}^{\ell}})^*
(  {R^{\ell-1}})^*
\end{equation*}
for $\ell=h,h-1,\dots,1$,
 starting from the middle matrix $M^h$ and
moving towards outer matrices.
Let
\begin{equation*}
{\bar{U}^L} = {U^L}  {C^L} ,
\end{equation*}
and
\begin{equation*}
 {\bar{V}^0}= {V^0}  {R^{0}},
\end{equation*}
then we have a further compressed factorization
\begin{equation*}
 K \approx  { \bar{U}^{L} \bar{G}^{L-1}\cdots \bar{G}^{h} \bar{M}^h  (\bar{H}^{h})^* \cdots (\bar{H}^{1})^*}  (\bar{V}^{0})^* ,
\end{equation*}
where all sparse factors have dense submatrices of size closer to
$r_0\times r_0$.
${\bar{U}^L} $ and $ {\bar{V}^0}$ are block-diagonal matrices and
the sizes of diagonal blocks depend on the distribution of
points in $X$ and $\Omega$. In the case of nonuniform distribution,
there might be diagonal blocks with size smaller than $r_0 \times r_0$
in ${\bar{U}^L} $ and $ {\bar{V}^0}$.
This motivates the sweep-in stage that contains another sequence of structure-preserving
matrix compression:
\begin{equation*}
{\bar{U}^L}\approx
{\dot{U}^L} \bar{C}^{L} ,
\end{equation*}
\begin{equation*}
 {\bar{V}^0} \approx {\dot{V}^0}
 \bar{R}^{0},
\end{equation*}
\begin{equation*}
\bar{C}^{\ell+1}  \bar{G}^\ell  \approx \dot{G}^\ell \bar{C}^{\ell}
\end{equation*}
for $\ell=L-1,L-2,\dots,h$,
\begin{equation*}
 ( \bar{H}^\ell)^* ( \bar{R}^{\ell-1})^* \approx
 ( \bar{R}^{\ell})^*  ( \dot{H}^\ell)^*
\end{equation*}
for $\ell=1,2,\dots,h$, starting from outer matrices to the middle matrix.
Finally, let $\dot{M}^h = \bar{C}^h \bar{M}^h(\bar{R}^h)^*$,
and one reaches the optimal IBF
\begin{equation}
\label{eqn:Kfactors}
 K \approx  { \dot{U}^{L} \dot{G}^{L-1}\cdots \dot{G}^{h} \dot{M}^h  (\dot{H}^{h})^* \cdots (\dot{H}^{1})^*}   (\dot{V}^{0})^*,
\end{equation}
where all sparse factors have dense submatrices of
nearly optimal size.

The optimal IBF represents $K$ as a product of $L+3$ sparse
matrices, where all factors are sparse matrices
with $\O(N)$ nonzero entries, and the prefactor of $O(N)$ is nearly optimal.
Once constructed,
the cost of applying $K$ to a given vector $g\in\bbC^N$ is $\O(N\log
N)$.

\subsection{Content}

The rest of this paper is organized as follows. Section \ref{sec:pre}
briefly reviews low-rank factorization techniques and
 the butterfly algorithm in \cite{FIO09}. Section
\ref{sec:preIBF} describes the one-dimensional preliminary
 interpolative butterfly factorization based on interpolative low-rank factorization.
 Section \ref{sec:optIBF} introduces a structure-preserving
 matrix compression technique and a sweeping method to
 further compress the preliminary IBF into an optimal one.
 Multidimensional extension to a general case when the phase function
 has singularity at $\xi=0$ is discussed in Section \ref{sec:mIBF}.
In Section \ref{sec:results}, numerical examples are
provided to demonstrate the efficiency of the proposed algorithms.
Finally, Section \ref{sec:conclusion} concludes this paper with a short discussion.

\section{Preliminary}
\label{sec:pre}

\subsection{Low-rank factorization}
\label{sec:LR}

This section reviews basic tools for efficient low-rank approximations
that are repeatedly used in this paper.

\textbf{Randomized low-rank approximation}

For a matrix $Z \in \bbC^{m\times n}$, we define a rank-$r$
approximate singular value decomposition (SVD) of $Z$ as
\[
Z \approx U_0 \Sigma_0 V_0^*,
\]
where $U_0\in \bbC^{m\times r}$ is unitary, $\Sigma_0\in \bbR^{r\times
  r}$ is diagonal, and $V_0\in \bbC^{n\times r}$ is unitary. A
straightforward method to obtain the optimal rank-$r$ approximation of
$Z$ is to compute its truncated SVD, where $U_0$ is the matrix with
the first $r$ left singular vectors, $\Sigma_0$ is a diagonal matrix
with the first $r$ singular values in decreasing order, and $V_0$ is
the matrix with the first $r$ right singular vectors.

The original truncated {SVD} of $Z$ takes
$\O(mn\min(m,n))$ operations. More efficient tools have been proposed
by introducing randomness in computing approximate SVDs for numerically
low-rank matrices. To name a few,  the one in \cite{Rec1} is
based on applying the matrix to random vectors while another one in
\cite{Butterfly4,symbol3} relies on sampling the matrix entries
randomly. Throughout this paper, the second one is applied to compute
large low-rank approximations because it only  takes linear operations
with respect to the matrix size. Readers are referred to \cite{Butterfly4,
symbol3} for detailed implementation.

When an approximate SVD $Z \approx U_0 \Sigma_0 V_0^*$ is ready, it
can be rearranged in several equivalent ways. First, one can write
\[
Z \approx U S V^*,
\]
where
\begin{equation}
  \label{eq:lowrankSVD}
  U=U_0\Sigma_0^{1/2},\, S = I \text{ and } V^*=\Sigma_0^{1/2}V_0^*,
\end{equation}
so that the left and right factors inherit similar singular
values of the original numerical low-rank matrix.
Depending on certain applications, sometimes it is better to write the
approximation as
\[
Z\approx UV^*
\]
where
\begin{equation}
  \label{eq:lowrankSVD2}
  U=U_0 \text{ and } V^*=\Sigma_0V_0^*,
\end{equation}
or
\begin{equation}
  \label{eq:lowrankSVD3}
  U=U_0 \Sigma_0\text{ and }
  V^*=V_0^*
\end{equation}
so that only one factor shares the singular values of $Z$.

\textbf{Interpolative low-rank approximation}

The randomized low-rank approximation is efficient if a kernel matrix is given, while
the interpolative low-rank approximation is efficient when the explicit formula
of the kernel function $K(x,\xi)=e^{2\pi \i \Phi(x,\xi)}$
is available, because the approximation can be constructed explicitly.

Following the theorems in \cite{FIO09,FIO14}, it can be proved that
the kernel function $K(x,\xi)=e^{2\pi \i \Phi(x,\xi)}$
satisfies the complementary low-rank property if it fulfills the properties
in Assumption \ref{asm}.
We will refresh the key idea here for one-dimensional
kernels. Let the set $X$ and $\Omega$
refer to the sets defined in \eqref{eqn:X} and \eqref{eqn:Omega}.
Let $A$ and $B$ be a box pair in the dyadic trees $T_X$ and $T_\Omega$
such that their levels satisfy $\ell_A + \ell_B =L$. Then a low-rank separated
representation
\[
K(x,\xi)=e^{2\pi \i \Phi(x,\xi)} \approx
\sum_{t=1}^{r_\eps} \alpha^{AB}_t(x) \beta^{AB}_t(\xi)\quad
\text{ for }x\in A, \xi\in B
\]
exists and can be constructed via the interpolative low-rank approximation as follows.

Let
\begin{equation}
  R^{AB}(x,\xi) := \Phi(x,\xi)-\Phi(c_A,\xi)-\Phi(x,c_B) +
  \Phi(c_A,c_B),
  \label{eqn:RAB}
\end{equation}
where $c_A$ and $c_B$ are the centers of $A$ and $B$
respectively, then the kernel can be written as
\begin{equation}
  e^{2\pi\i \Phi(x,\xi)} =
  e^{2\pi\i \Phi(c_A,\xi)} e^{2\pi\i \Phi(x,c_B)}
  e^{-2\pi\i\Phi(c_A,c_B)} e^{2\pi\i R^{AB}(x,\xi)}.
  \label{eqn:terms}
\end{equation}
Hence, the low-rank approximation of $K(x,\xi)$ in $A\times B$ is reduced to
the low-rank approximation of $e^{2\pi\i R^{AB}(x,\xi)}$.

\newcommand{\talpha}{\tilde{\alpha}}
\newcommand{\tbeta}{\tilde{\beta}}

Let $w_A$ and $w_B$ denote the lengths of
intervals $A$ and $B$, respectively. A Lagrange interpolation
with Chebyshev points in $x$ when $w_A \le 1/\sqrt{N}$ and in
$\xi$ when $w_B \le \sqrt{N}$ is applied to construct the
low-rank approximation of $e^{2\pi\i R^{AB}(x,\xi)}$.
For this purpose, we associate with
each interval a Chebyshev grid as follows.

For a fixed integer $r$, the
Chebyshev grid of order $r$ on
$[-1/2,1/2]$ is defined by
\[
\left\{ z_t = \frac{1}{2} \cos \left( \frac{t\pi}{r-1} \right) \right\}_{0\le t\le r-1}.
\]
A grid {\em adapted to} an interval $A$ with center $c_A$ and
length $w_A$ is then defined via shifting and scaling as
\[
\{x_t\}_{t=0,1,\ldots, r-1}=\{c_A + w_A z_t\}_{t=0,1,\ldots, r-1}.
\]
Given a set of grid points $\{x_t\}_{t=0,1,\ldots, r-1}$ in $A$, define
Lagrange interpolation polynomials $M^A_t(x)$ taking value $1$
at $x_t$ and $0$ at the other Chebyshev grid points
\[
M^A_t(x)=\prod_{0\leq j\leq r-1,j\neq t} \frac{x-x_j}{x_t-x_j}.
\]
Similarly, $M_t^B$ is denoted as the Lagrange interpolation
polynomials for the interval $B$.

Now we are ready to construct the low-rank approximation of $e^{2\pi\i R^{AB}(x,\xi)}$ with $r_\eps$ Chebyshev points for $\eps$-accuracy by interpolation:
\begin{itemize}
  \item when $w_B \le \sqrt{N}$, the Lagrange interpolation of
    $e^{2\pi\i R^{AB}(x,\xi)}$ in $\xi$ on a Chebyshev grid
     $\{g^B_t\}_{1\le t \le r_\eps}$ adapted to $B$
    obeys
    \begin{equation}
   e^{2\pi\i R^{AB}(x,\xi)} \approx \sum_{t=1}^{r_\eps} e^{2\pi\i R^{AB}(x,g^B_t)} M^B_t(\xi) , \quad \forall x\in A, \forall \xi\in B,
      \label{eqn:intp1}
    \end{equation}
  \item when $w_A \le 1/\sqrt{N}$, the Lagrange interpolation of
    $e^{2\pi\i R^{AB}(x,\xi)}$ in $x$ on a
    Chebyshev grid $\{g^A_t\}_{1\le t\le r_\eps}$ adapted to $A$ obeys
    \begin{equation}
   e^{2\pi\i R^{AB}(x,\xi)} \approx \sum_{t=1}^{r_\eps} M^A_t(x) e^{2\pi\i R^{AB}(g^A_t,\xi)} , \quad \forall x\in A, \forall \xi\in B.
      \label{eqn:intp2}
    \end{equation}
  \end{itemize}

Finally, we are ready to construct the low-rank approximation for the
kernel $e^{2\pi\i \Phi(x,\xi)}$:
\begin{itemize}
\item when $w_B \le \sqrt{N}$, we multiply \eqref{eqn:intp1} with
$e^{2\pi\i \Phi(c_A,\xi)} e^{2\pi\i \Phi(x,c_B)}
e^{-2\pi\i\Phi(c_A,c_B)}$, which gives that $\forall x\in A, \forall
\xi\in B$
\begin{equation}
\label{eqn:ab1}
e^{2\pi\i\Phi(x,\xi)} \approx \sum_{t=1}^{r_\eps} e^{2\pi\i\Phi(x,g^B_t)}
\left(
e^{-2\pi\i\Phi(c_A,g^B_t)} M^B_t(\xi) e^{2\pi\i\Phi(c_A,\xi)}
\right);
\end{equation}
\item when $w_A \le
1/\sqrt{N}$, multiply \eqref{eqn:intp2} with $e^{2\pi\i \Phi(c_A,\xi)}
e^{2\pi\i\Phi(x,c_B)} e^{-2\pi\i\Phi(c_A,c_B)}$ and obtain that
$\forall x\in A, \forall \xi\in B$
\begin{equation}
\label{eqn:ab2}
  e^{2\pi\i\Phi(x,\xi)} \approx \sum_{t=1}^{r_\eps}
  \left(
  e^{2\pi\i\Phi(x,c_B)} M^A_t(x) e^{-2\pi\i\Phi(g^A_t,c_B)}
  \right)
  e^{2\pi\i\Phi(g^A_t,\xi)}.
\end{equation}
\end{itemize}

The interpolative low-rank factorization can be constructed on-the-fly
from the explicit formulas above, which is the main advantage over
randomized low-rank approximations. However, since it relies on
the information on a fixed Chebyshev grid, the number of Chebyshev
points must be sufficiently large to obtain an accurate approximation,
i.e., the {\em $\eps$-separation rank} $r_\eps$
might be greater than the true numerical rank with $\eps$ accuracy.

\subsection{Butterfly algorithm}
\label{sec:BA}

This section provides a brief description of the overall structure of
the butterfly algorithm based on the interpolative low-rank approximation
in the previous section.
In this section, $X$ and $\Omega$ refer to
two general sets of $N$ points in $\R$, respectively. With no loss of
generality, we assume the
points in these two sets are distributed quasi-uniformly but they are
not necessarily the sets defined in \eqref{eqn:X} and
\eqref{eqn:Omega}.

Given an input $\{g(\xi), \xi \in \Omega\}$, the goal is to compute
the potentials $\{u(x), x\in X\}$ defined by
\[
u(x) = \sum_{\xi\in \Omega} K(x,\xi) g(\xi), \quad x \in X,
\]
where $K(x,\xi)$ is a kernel function. The main data structure of the butterfly
algorithm is a pair of dyadic trees $T_X$ and $T_\Omega$. Recall that
for any pair of intervals $A\times B \in T_X\times T_\Omega$
obeying the condition $\ell_A + \ell_B = L$,
 the submatrix $\{K(x,\xi)\}_{x \in A, \xi \in B}$ is
approximately of a constant rank. An explicit method to construct
its low-rank approximation is given by the interpolative low-rank approximation.
More precisely, for any $\eps>0$,
there exists a constant $r_\eps$ independent of $N$ and two sets of
functions $\{\alpha^{AB}_t (x)\}_{1\le t \le r_\eps}$ and
$\{\beta^{AB}_t (\xi)\}_{1\le t \le r_\eps}$  given in
\eqref{eqn:ab1} or \eqref{eqn:ab2} such that
\begin{equation}
  \left| K(x,\xi) - \sum_{t=1}^{r_\eps} \alpha^{AB}_t(x)
  \beta^{AB}_t(\xi) \right| \le \eps, \quad \forall x\in A, \forall
  \xi\in B.
  \label{eqn:glr}
\end{equation}

\begin{figure}
\centering
\begin{tikzpicture}[scale=0.3]

\draw [->,double,double distance=1pt] (-6,9) -> (-6,1);

\draw (0,10) -- (-4,0) -- (4,0) -- cycle;

\draw [->,double,double distance=1pt] (20,1) -> (20,9);

\draw (14,10) -- (10,0) -- (18,0) -- cycle;

\draw [latex-latex,line width=1pt] (9,5) -> (5,5);
\draw [latex-latex,line width=1pt] (9,9) -> (5,1);
\draw [latex-latex,line width=1pt] (9,1) -> (5,9);

\coordinate [label=above:$T_X$] (X) at (0,10);
\coordinate [label=above:$T_\Omega$] (X) at (14,10);

\draw (-2,5) -- (2,5);
\draw (12,5) -- (16,5);

\coordinate [label=right:$\frac{L}{2}$] (X) at (16,5);
\coordinate [label=left:$\frac{L}{2}$] (X) at (-2,5);
\end{tikzpicture}
\caption{Trees of the row and column indices.
    Left: $T_X$ for the row indices $X$.
    Right: $T_\Omega$ for the column indices $\Omega$.
    The interaction between $A\in T_X$ and $B\in T_\Omega$
    starts at the root of $T_X$ and the leaves of $T_\Omega$. }
\label{fig:domain-tree-BA}
\end{figure}
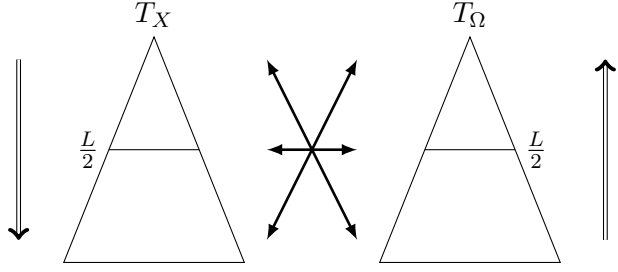

For a given interval $B$ in $\Omega$, define $u^B(x)$ to be the {\em
  restricted potential} over the sources $\xi\in B$
\[
u^B(x) = \sum_{\xi\in B} K(x,\xi) g(\xi).
\]
The low-rank property gives a compact expansion for $\{u^B(x)\}_{x\in
  A}$ as summing \eqref{eqn:glr} over $\xi\in B$ with coefficients $g(\xi)$
gives
\[
\left| u^B(x) - \sum_{t=1}^{r_\eps} \alpha^{AB}_t(x) \left( \sum_{\xi\in B} \beta^{AB}_t(\xi) g(\xi) \right) \right|
\le \left( \sum_{\xi\in B} |g(\xi)| \right) \eps,
\quad \forall x \in A.
\]
Therefore, if one can find coefficients $\{\coef{AB}{t}\}_{1\le t\le
  r_\eps}$ obeying
\begin{equation}
  \coef{AB}{t} \approx \sum_{\xi\in B} \beta^{AB}_t(\xi) g(\xi), \quad 1\le t\le r_\eps,
  \label{eqn:delta}
\end{equation}
then the restricted potential $\{u^B(x)\}_{x\in A}$ admits a compact
expansion
\[
\left| u^B(x) - \sum_{t=1}^{r_\eps} \alpha^{AB}_t(x) \coef{AB}{t} \right| \le \left( \sum_{\xi\in B} |g(\xi)| \right) \eps,
\quad \forall x\in A.
\]
The butterfly algorithm below provides an efficient way for
computing $\{\coef{AB}{t}\}_{1\le t\le r_\eps}$ recursively.  The
general structure of the algorithm consists of a top-down traversal of
$T_X$ and a bottom-up traversal of $T_\Omega$, carried out
simultaneously. A schematic illustration of the data flow
in this algorithm is provided in Figure
\ref{fig:domain-tree-BA}.
\begin{algo}{Butterfly algorithm}
\label{aglo:BA}

\begin{enumerate}
\item{\emph{Preliminaries.}} Construct the trees $T_X$ and $T_\Omega$.

\item{\emph{Initialization.}} Let $A$ be the root of $T_X$. For each leaf interval $B$ of
  $T_\Omega$, construct the expansion coefficients $\{
  \coef{AB}{t}\}_{1\le t \le r_\eps}$ for the potential
  $\{u^B(x)\}_{x\in A}$ by simply setting
  \begin{equation}
    \coef{AB}{t} = \sum_{\xi\in B} \beta^{AB}_t(\xi) g(\xi), \quad 1\le t \le r_\eps.
  \label{eqn:bf1}
  \end{equation}
    By the interpolative low-rank approximation, we can define the expansion coefficients $\{\coef{AB}{t}\}_{1\le
    t\le r_\eps}$ by
  \begin{equation}
    \label{eqn:1}
    \coef{AB}{t} :=
    e^{-2\pi\i \Phi(c_A,g_t^B)} \sum_{\xi\in B}\left( M_t^B(\xi)
    e^{2\pi\i\Phi(c_A,\xi)}g(\xi)\right).
  \end{equation}

  \item{\emph{Recursion.}} For $\ell = 1, 2, \ldots, L/2$, visit level $\ell$ in $T_X$ and
  level $L-\ell$ in $T_\Omega$. For each pair $(A,B)$ with $\ell_A =
  \ell$ and $\ell_B = L-\ell$, construct the expansion coefficients
  $\{\coef{AB}{t}\}_{1\le t \le r_\eps}$ for the potential
  $\{u^B(x)\}_{x\in A}$ using the low-rank representation constructed
  at the previous level. Let
  $P$ be $A$'s parent and $C$ be a child of $B$. Throughout, we shall
  use the notation $C\succ B$ when $C$ is a child of $B$. At level
  $\ell-1$, the expansion coefficients $\{\coef{PC}{s}\}_{1\le
    s\le r_\eps}$ of $\{u^{C}(x)\}_{x\in P}$ are readily available
  and we have
  \[
  \left| u^{C}(x) - \sum_{s=1}^{r_\eps} \alpha^{PC}_{s}(x) \coef{PC}{s} \right| \le \left( \sum_{\xi\in C} |g(\xi)| \right) \eps,
  \quad \forall x\in P.
  \]
  Since $u^B(x) = \sum_{C\succ B} u^{C}(x)$, the previous inequality
  implies that
  \[
  \left| u^B(x) - \sum_{C\succ B} \sum_{s=1}^{r_\eps} \alpha^{PC}_{s}(x) \coef{PC}{s} \right| \le \left( \sum_{\xi\in B} |g(\xi)| \right) \eps,
  \quad \forall x\in P.
  \]
  Since $A \subset P$, the above approximation is of course true for
  any $x \in A$. However, since $\ell_A + \ell_B = L$, the sequence of
  restricted potentials $\{u^B(x)\}_{x\in A}$ also has a low-rank
  approximation of size $r_\eps$, namely,
  \[
  \left| u^B(x) - \sum_{t=1}^{r_\eps} \alpha^{AB}_t(x) \coef{AB}{t} \right| \le \left( \sum_{\xi\in B} |g(\xi)| \right) \eps,
  \quad \forall x\in A.
  \]
  Combining the last two approximations, we obtain that
  $\{\coef{AB}{t}\}_{1\le t\le r_\eps}$ should obey
  \begin{equation}
    \sum_{t=1}^{r_\eps} \alpha^{AB}_t(x) \coef{AB}{t} \approx
    \sum_{C\succ B} \sum_{s=1}^{r_\eps} \alpha^{PC}_{s}(x) \coef{PC}{s}, \quad \forall x\in A.
    \label{eqn:bf2}
  \end{equation}
  This is an over-determined linear system for
  $\{\coef{AB}{t}\}_{1\le t\le r_\eps}$ when
  $\{\coef{PC}{s}\}_{1\le s\le r_\eps,C\succ B}$ are available.
  The butterfly algorithm uses an
  efficient linear transformation approximately mapping
  $\{\coef{PC}{s}\}_{1\le s\le r_\eps,C\succ B}$ into
  $\{\coef{AB}{t}\}_{1\le t\le r_\eps}$ as follows
  \begin{equation}
    \label{eqn:2}
    \coef{AB}{t} := e^{-2\pi\i\Phi(c_A,g_t^B)}\sum_{C\succ B}\sum_{s=1}^{r_\eps} M_t^B(g_{s}^{C})
    e^{2\pi\i\Phi(c_A,g_{s}^{C})}\coef{PC}{s}.
  \end{equation}

\item{\emph{Switch.}} For the levels visited, the Chebyshev
  interpolation is applied in variable $\xi$, while the interpolation
  is applied in variable $x$ for levels $\ell>L/2$.  Hence, we
  are switching the interpolation method at this step.  Now we are
  still working on level $\ell=L/2$ and the same domain pairs
  $(A,B)$ in the last step.  Let $\coef{AB}{s}$ denote the expansion
  coefficients obtained by Chebyshev interpolation in variable $\xi$ in the
  last step.  Correspondingly, $\{g_s^B\}_s$ are the grid points in
  $B$ in the last step.  We take advantage of the interpolation in
  variable $x$ in $A$ and generate grid points $\{g_t^A\}_{1\le t\le
    r_\eps}$ in $A$.  Then we can define new expansion coefficients
  \[
  \coef{AB}{t} := \sum_{s=1}^{r_\eps} e^{2\pi\i\Phi(g_t^A,g_s^B)}\coef{AB}{s}.
  \]
\item{\emph{Recursion.}} Similar to the discussion in Step $3$,
  we go up in tree $T_{\Omega}$ and down in tree
  $T_X$ at the same time until we reach the level
  $\ell=L$. We construct the approximation functions by Chebyshev
  interpolation in variable $x$ as follows:
  \begin{eqnarray}
    \label{eqn:intx}
    \alpha_t^{AB}(x)= e^{2\pi\i\Phi(x,c_B)}M_t^A(x) e^{-2\pi\i\Phi(g_t^A,c_B)},
    &\beta_t^{AB}(\xi)=e^{2\pi\i\Phi(g_t^A,\xi)}.
  \end{eqnarray}
  Hence, the new expansion coefficients $\{\coef{AB}{t}\}_{1\le
    t\le r_\eps}$ can be defined as
  \begin{equation}
    \label{eqn:3}
    \coef{AB}{t} := \sum_{C\succ B} e^{2\pi\i\Phi(g_t^A,c_{C})}
    \sum_{s=1}^{r_\eps}\left(M_{s}^{P}(g_{t}^{A})
    e^{-2\pi\i\Phi(g_{s}^{P},c_{C})}\coef{PC}{s}\right),
  \end{equation}
  where again $P$ is $A$'s parent and $C$ is a child interval of $B$.
\item{\emph{Termination.}} Finally, $\ell = L$ and set $B$ to be the root node of
  $T_\Omega$. For each leaf interval $A \in T_X$, use the constructed
  expansion coefficients $\{\coef{AB}{t}\}_{1\le t\le r_\eps}$  in \eqref{eqn:3} to
  evaluate $u^B(x)$ for each $x \in A$,
  \begin{equation}
    \begin{split}
    u(x) = u^B(x)&= \sum_{t=1}^{r_\eps} \alpha^{AB}_t (x) \coef{AB}{t}\\
    &=e^{2\pi\i\Phi(x,c_B)} \sum_{t=1}^{r_\eps}\left(M_t^A(x)e^{-2\pi\i\Phi(g_t^A,c_B)} \coef{AB}{t} \right).
    \label{eqn:bf3}
    \end{split}
  \end{equation}

\end{enumerate}

\end{algo}

\section{Preliminary interpolative  butterfly factorization (IBF)}
\label{sec:preIBF}

This section presents the preliminary interpolative butterfly factorization (IBF) for a
matrix $K\in\bbC^{N\times N}$ when its kernel function satisfies
Assumption \ref{asm}. In fact, the preliminary interpolative butterfly factorization is 
a matrix representation of the butterfly algorithm \cite{FIO09}. 
Similar to the butterfly algorithm, we adopt the same
notation of point sets $X$ and $\Omega$, trees $T_X$ and $T_\Omega$ of depth
$L$ (assumed to be an even number). At each level
$\ell$, $\ell=0,\dots,L$, we denote the $i$th node at
level $\ell$ in $T_X$ as $A^\ell_i$ for $i=0,1,\dots,2^\ell-1$ and the
$j$th node at level $L-\ell$ in $T_\Omega$ as $B^{L-\ell}_j$ for
$j=0,1,\dots,2^{L-\ell}-1$. These nodes naturally partition $K$ into
$\O(N)$ submatrices $K_{A^\ell_i,B^{L-\ell}_j}$. For simplicity, we
write $K^\ell_{i,j}:=K_{A^\ell_i,B^{L-\ell}_j}$, where the superscript
is used to indicate the level (in $T_X$). With abuse of notation,
sometimes we use the subscripts ${i,j}$ and the superscript $\ell$ on a matrix or a vector
corresponding the domain pair $(A^\ell_i,B^{L-\ell}_j)$ for simplicity,
although it is not a submatrix or a subvector restricted in
$(A^\ell_i,B^{L-\ell}_j)$.

The preliminary IBF is based on the observation that
the butterfly algorithm in Section \ref{sec:BA} can be written in a form of matrix
factorization. The operations in Step 2 to 6 in Algorithm \ref{aglo:BA} are
essentially a sequence of matrix vector multiplications with $O(\log N)$
multiplications, each matrix of which has only $O(r_\eps^2 N)$ nonzero
entries. Since all the operations are given explicitly based on Chebyshev
interpolation, these sparse matrices can be formed by explicit formulas. By
formulating these matrices step by step following the flow in Algorithm
\ref{aglo:BA}, the preliminary IBF is proposed
as follows.

\begin{enumerate}
\item{\emph{Preliminaries.}} Construct the trees $T_X$ and $T_\Omega$.

\item{\emph{Initialization.}} At level $\ell=0$, for each $j=0,1,\dots,2^{L}-1$, let $A$ be $\Omega$ and $B$ be
$B^{L}_j$ of  $T_\Omega$. Construct the expansion coefficients $\{
  \coef{AB}{t}\}_{1\le t \le r_\eps}$ for the potential
  $\{u^B(x)\}_{x\in A}$ by simply setting
  \begin{equation}
    \label{eqn:PIBF1}
    \coef{AB}{t} :=
    e^{-2\pi\i \Phi(c_A,g_t^B)} \sum_{\xi\in B}\left( M_t^B(\xi)
    e^{2\pi\i\Phi(c_A,\xi)}g(\xi)\right).
  \end{equation}

  For each  $j=0,1,\dots,2^{L}-1$, a column vector $\Coef{0}{j}$ corresponding
  to the domain pair $(A,B)=(\Omega,B^{L}_j)$ is defined
  as
\begin{equation*}
\Coef{0}{j} =
\begin{pmatrix}
  \coef{AB}{1}\\
  \vdots\\
  \coef{AB}{r_\eps}
  \end{pmatrix}.
\end{equation*}
Let $(V^0_j)^*\in \mathbb{C}^{r_\eps\times O(1)}$ represent the linear transformation in \eqref{eqn:PIBF1}
and $g^0_j$ be the vector representing $g(\xi)$ for $\xi\in B^L_j$. Then we have
\begin{equation}\label{eqn:step2mv1}
 \Coef{0}{j}  =(V^0_{j})^*  g^0_{j}.
\end{equation}
  As we shall see later, the conjugate transpose $^*$ is
  applied for the purpose of notation consistency.
  By assembling the matrix-vector multiplications in \eqref{eqn:step2mv1},
  all the operations in this step can be written as
   \[
   \Coef{0}{}=
   \begin{pmatrix}
   \Coef{0}{0} \\
   \vdots \\
   \Coef{0}{2^L-1}
  \end{pmatrix}
  =(V^0)^*g,
   \]
   where
   \[V^0=\diag\left\{V^0_{0},V^0_{1},
 \cdots,V^0_{2^{L}-1}\right\}.\]

  \item{\emph{Recursion.}} For $\ell = 1, 2, \ldots, L/2$,
  visit level $\ell$ in $T_X$ and
  level $L-\ell$ in $T_\Omega$.
  At each level $\ell$, for each pair $(A,B)$ with $\ell_A =
  \ell$ and $\ell_B = L-\ell$, construct the expansion coefficients
  $\{\coef{AB}{t}\}_{1\le t \le r_\eps}$
  using the low-rank representation constructed
  at the previous level ($\ell = 0$ is the initialization step). Let
  $P$ be $A$'s parent and $C$ be a child of $B$. An
  efficient linear transformation approximately mapping
  $\{\coef{PC}{s}\}_{1\le s\le r_\eps,C\succ B}$ into
  $\{\coef{AB}{t}\}_{1\le t\le r_\eps}$ is constructed by
  \begin{equation}
    \label{eqn:PIBF2}
    \coef{AB}{t} := e^{-2\pi\i\Phi(c_A,g_t^B)}\sum_{C\succ B}\sum_{s=1}^{r_\eps} M_t^B(g_{s}^{C})
    e^{2\pi\i\Phi(c_A,g_{s}^{C})}\coef{PC}{s}.
  \end{equation}

  At each level $\ell$, for each $i=0,1,\dots,2^\ell-1$
  and $j=0,1,\dots,2^{L-\ell}-1$,
  a column vector $\Coef{\ell}{i,j}$ is defined as
\begin{equation}
\Coef{\ell}{i,j}=
\begin{pmatrix}
  \coef{AB}{1}\\
  \vdots\\
  \coef{AB}{r_\eps}
  \end{pmatrix} \in \mathbb{C}^{r_\eps},
\end{equation}
  where the domain pair $(A,B)=(A^\ell_i,B^{L-\ell}_j)$.
  For a fixed $j$,  stacking vectors $\{\Coef{\ell}{i,j}\}_i$ together
  forms a larger column vector $\Coef{\ell}{j}$
  and stacking vectors $\{\Coef{\ell}{j}\}_j$ together
  forms a column vector $\Coef{\ell}{}$,
  \begin{equation*}
  \Coef{\ell}{j} =
  \begin{pmatrix}
  \Coef{\ell}{0,j}\\
  \vdots \\
  \Coef{\ell}{2^\ell-1,j}
  \end{pmatrix},
  \quad
  \Coef{\ell}{} =
  \begin{pmatrix}
  \Coef{\ell}{0}\\
  \vdots \\
  \Coef{\ell}{2^{L-\ell}-1}
  \end{pmatrix}.
\end{equation*}

The linear transformation in \eqref{eqn:PIBF2} maps two small column vectors
$\Coef{\ell-1}{i,2j}$ and $\Coef{\ell-1}{i,2j+1}$ at the previous level $\ell-1$
to a small column vector $\Coef{\ell}{2i,j}$
if $A$ is the first child of $P$
(or $\Coef{\ell}{2i+1,j}$ if $A$ is the second child of $P$)
at the current level $\ell$ for all $i=0,1,\dots,2^{\ell-1}-1$
and $j=0,1,\dots,2^{L-\ell}-1$.
Hence, for each pair $(i,j)$, the linear transformation in \eqref{eqn:PIBF2}
can be written as a matrix vector multiplication
   \begin{equation}
  \label{eqn:PIBFMatR1}
  \Coef{\ell}{2i,j}=
  \begin{pmatrix}
  H^{\ell}_{2i,2j} &  H^{\ell}_{2i,2j+1}
  \end{pmatrix}
   \begin{pmatrix}
 \Coef{\ell-1}{i,2j}\\
  \Coef{\ell-1}{i,2j+1}
  \end{pmatrix},
  \end{equation}
  where $H^{\ell}_{2i,2j}\in \mathbb{C}^{r_\eps\times r_\eps}$
  and $H^{\ell}_{2i,2j+1}\in \mathbb{C}^{r_\eps\times r_\eps}$, or
   \begin{equation}
  \label{eqn:PIBFMatR1}
  \Coef{\ell}{2i+1,j}=
  \begin{pmatrix}
  H^{\ell}_{2i+1,2j} & H^{\ell}_{2i+1,2j+1}
  \end{pmatrix}
   \begin{pmatrix}
  \Coef{\ell-1}{i,2j}\\
  \Coef{\ell-1}{i,2j+1}
  \end{pmatrix},
  \end{equation}
    where $H^{\ell}_{2i+1,2j}\in \mathbb{C}^{r_\eps\times r_\eps}$
    and $H^{\ell}_{2i+1,2j+1}\in \mathbb{C}^{r_\eps\times r_\eps}$.
Assemble the above small matrices together to get
\begin{equation*}
H^{\ell}_{j}=
\begin{pmatrix}
H^{\ell}_{0,2j} & H^{\ell}_{1,2j} & & &\\
 & & H^{\ell}_{2,2j} & H^{\ell}_{3,2j} & &\\
 & & & & \ddots & &\\
 & & & & & H^{\ell}_{2^{\ell}-2,2j} & H^{\ell}_{2^{\ell}-1,2j}\\
\midrule
H^{\ell}_{0,2j+1} & H^{\ell}_{1,2j+1} & & &\\
 & & H^{\ell}_{2,2j+1} & H^{\ell}_{3,2j+1} & &\\
 & & & & \ddots &\\
 & & & & & H^{\ell}_{2^{\ell}-2,2j+1} & H^{\ell}_{2^{\ell}-1,2j+1}
\end{pmatrix}
\end{equation*}
for $j=0,1,\dots,2^{L-\ell}-1$ and
  \begin{equation*}
H^{\ell}=
\begin{pmatrix}
H^{\ell}_{0} & & &\\
 & H^{\ell}_{1} & &\\
 & & \ddots &\\
 & & & H^{\ell}_{2^{L-\ell}-1}\\
\end{pmatrix},
\end{equation*}
then all the operations in this step can be represented by a large matrix-vector multiplication
 \[
 \Coef{\ell}{}=(H^\ell)^* \Coef{\ell-1}{}
 \]
  for $\ell = 1, 2, \ldots, L/2$.

\item{\emph{Switch.}} For the levels visited, the Chebyshev
  interpolation is applied in variable $\xi$, while the interpolation
  is applied in variable $x$ for levels $\ell>L/2$.  Hence, we
  are switching the interpolation method at this step.  Now we are
  still working at level $\ell=L/2$ and the same domain pairs
  $(A,B)=(A^{L/2}_{i},B^{L/2}_j)$ in the last step.  Recall that
  $\Coef{L/2}{}$ denote the expansion
  coefficients obtained by Chebyshev interpolation in variable $\xi$ in the
  last step.  Correspondingly, $\{g_s^B\}_s$ are the grid points in
  $B$ in the last step.  We take advantage of the interpolation in
  variable $x$ in $A$ and generate grid points $\{g_t^A\}_{1\le t\le
    r_\eps}$ in $A$.  Then we can define new expansion coefficients
    for Chebyshev interpolation in variable $x$ for future steps as
  \begin{equation}
  \label{eqn:switchlocM}
  \coef{AB}{t} := \sum_{s=1}^{r_\eps} e^{2\pi\i\Phi(g_t^A,g_s^B)}\coef{AB}{s}.
  \end{equation}
 Let $\tilde{M}^{L/2}_{i,j}\in \mathbb{C}^{r_\eps\times r_\eps}$ represent the linear transformation introduced by $\{ e^{2\pi\i\Phi(g_t^A,g_s^B)} \}_{1\leq t\leq r_\eps,1\leq s\leq r_\eps}$ for $i$, $j=0,1,\cdots,2^{L/2}-1$. Then \eqref{eqn:switchlocM} is equivalent to
 \[
 \tCoef{L/2}{i,j}=\tilde{M}^{L/2}_{i,j} \Coef{L/2}{i,j}
 \]
 for each domain pair $(A,B)=(A^{L/2}_{i},B^{L/2}_j)$. Recall that
  \begin{equation*}
 \Coef{L/2}{j} =
  \begin{pmatrix}
  \Coef{L/2}{0,j} \\
  \vdots \\
  \Coef{L/2}{2^{L/2}-1,j}
  \end{pmatrix},
  \text{ and }
  \Coef{L/2}{} =
  \begin{pmatrix}
  \Coef{L/2}{0} \\
  \vdots \\
  \Coef{L/2}{2^{L/2}-1}
  \end{pmatrix}.
\end{equation*}
If we define the expansion coefficients in a new order by
  \begin{equation*}
  \tCoef{L/2}{i}=
  \begin{pmatrix}
  \tCoef{L/2}{i,0}\\
  \vdots \\
  \tCoef{L/2}{i,2^{L/2}-1}
  \end{pmatrix},
  \text{ and }
  \tCoef{L/2}{}=
  \begin{pmatrix}
  \tCoef{L/2}{0}\\
  \vdots \\
  \tCoef{L/2}{2^{L/2}-1}
  \end{pmatrix},
\end{equation*}
then all the operations in \eqref{eqn:switchlocM} for all domain pairs can be represented by a large matrix-vector multiplication
\[
\tCoef{L/2}{}=M^{L/2}\Coef{L/2}{},
\]
where
\[
M^{L/2}=
   \begin{pmatrix}
      M^{L/2}_{0,0} & M^{L/2}_{0,1} & \cdots &M^{L/2}_{0,2^{L/2}-1}\\
     M^{L/2}_{1,0} &  M^{L/2}_{1,1} & & M^{L/2}_{1,2^{L/2}-1}\\
      \vdots & & \ddots & \\
     M^{L/2}_{2^{L/2}-1,0} & M^{L/2}_{2^{L/2}-1,1} & & M^{L/2}_{2^{L/2}-1,2^{L/2}-1}
    \end{pmatrix}\in \mathbb{C}^{r_\eps N\times r_\eps N}
\]
and $M^{L/2}_{i,j}\in \mathbb{C}^{r_\eps \sqrt{N}\times r_\eps \sqrt{N}}$ is also a $2^{L/2}\times 2^{L/2}$ block matrix with the only nonzero block at the location $(j,i)$ being $\tilde{M}^{L/2}_{i,j}\in \mathbb{C}^{r_\eps\times r_\eps}$.
By the abuse of notation, we drop the tilde notation $\tilde{\cdot}$ of the expansion coefficients
in the second half of the algorithm description.

\item{\emph{Recursion.}} Similar to the discussion in Step $3$,
  we go up in tree $T_{\Omega}$ and down in tree
  $T_X$ at the same time for all level
  $\ell=L/2+1,\cdots,L$. At each level $\ell$,
   for any domain pair $(A,B)=(A^\ell_{i},B^{L-\ell}_j)$,
   for $i=0,1,\dots,2^{\ell}-1$ and $j=0,1,\dots,2^{L-\ell}-1$,
   we construct the new expansion coefficients
   $\{\coef{AB}{t}\}_{1\le t\le r_\eps}$ by
  \begin{equation}
    \label{eqn:PIBF3}
    \coef{AB}{t} := \sum_{C\succ B} e^{2\pi\i\Phi(g_t^A,c_{C})}
    \sum_{s=1}^{r_\eps}\left(M_{s}^{P}(g_{t}^{A})
    e^{-2\pi\i\Phi(g_{s}^{P},c_{C})}\coef{PC}{s}\right),
  \end{equation}
  where again $P$ is $A$'s parent and $C$ is a child interval of $B$.
  The coefficients are assembled as
  \begin{equation*}
  \Coef{\ell}{i,j} =
  \begin{pmatrix}
  \coef{AB}{1}\\
  \vdots\\
  \coef{AB}{r_\eps}
  \end{pmatrix},
  \quad
  \Coef{\ell}{i} =
  \begin{pmatrix}
  \Coef{\ell}{i,0}\\
  \vdots\\
  \Coef{\ell}{i,2^{L-\ell}-1}
  \end{pmatrix},
  \text{ and }
  \Coef{\ell}{} =
  \begin{pmatrix}
  \Coef{\ell}{0}\\
  \vdots\\
  \Coef{\ell}{2^{\ell}-1}
  \end{pmatrix},
  \end{equation*}
for $i=0,1,\dots,2^{\ell}-1$, $j=0,1,\dots,2^{L-\ell}-1$.

Similarly to the first recursion step,
the linear transformation in \eqref{eqn:PIBF3} maps
two small column vectors
$\Coef{\ell-1}{i,2j}$ and $\Coef{\ell-1}{i,2j+1}$
at the previous level $\ell-1$
to a small column vector $\Coef{\ell}{2i,j}$
if $A$ is the first child of $P$
(or $\Coef{\ell}{2i+1,j}$ if $A$ is the second child of $P$)
at the current level $\ell$
for all $i=0,1,\dots,2^{\ell-1}-1$ and $j=0,1,\dots,2^{L-\ell}-1$.
Hence, for each pair $(i,j)$,
the linear transformation in \eqref{eqn:PIBF3}
can be written as a matrix vector multiplication
  \begin{equation}
  \label{eqn:PIBFMatR1}
  \Coef{\ell}{2i,j}=
  \begin{pmatrix}
 G^{\ell-1}_{2i,2j} &  G^{\ell-1}_{2i,2j+1}
  \end{pmatrix}
   \begin{pmatrix}
 \Coef{\ell-1}{i,2j}\\
  \Coef{\ell-1}{i,2j+1}
  \end{pmatrix},
  \end{equation}
  where $G^{\ell-1}_{2i,2j}\in \mathbb{C}^{r_\eps\times r_\eps}$
  and $G^{\ell-1}_{2i,2j+1}\in \mathbb{C}^{r_\eps\times r_\eps}$, or
   \begin{equation}
  \label{eqn:PIBFMatR1}
  \Coef{\ell}{2i+1,j}=
  \begin{pmatrix}
 G^{\ell-1}_{2i+1,2j} & G^{\ell-1}_{2i+1,2j+1}
  \end{pmatrix}
   \begin{pmatrix}
 \Coef{\ell-1}{i,2j}\\
  \Coef{\ell-1}{i,2j+1}
  \end{pmatrix},
  \end{equation}
where $ G^{\ell-1}_{2i+1,2j}\in \mathbb{C}^{r_\eps\times r_\eps}$
and $ G^{\ell-1}_{2i+1,2j+1}\in \mathbb{C}^{r_\eps\times r_\eps}$.
    Assemble the above small matrices together to get
    \begin{equation*}
G^{\ell-1}_{i}=
\begin{pmatrix}
G^{\ell-1}_{2i,0} & G^{\ell-1}_{2i,1} & & &\\
 & & G^{\ell-1}_{2i,2} & G^{\ell-1}_{2i,3} & &\\
 & & & & \ddots &\\
 & & & & & G^{\ell-1}_{2i,2^{L-\ell}-2} & G^{\ell-1}_{2i,2^{L-\ell}-1}\\
\midrule
G^{\ell-1}_{2i+1,0} & G^{\ell-1}_{2i+1,1} & & &\\
 & & G^{\ell-1}_{2i+1,2} & G^{\ell-1}_{2i+1,3} & &\\
 & & & & \ddots &\\
 & & & & & G^{\ell-1}_{2i+1,2^{L-\ell}-2}
  & G^{\ell-1}_{2i+1,2^{L-\ell}-1}
\end{pmatrix}
\end{equation*}
for $i=0,1,\dots,2^{\ell-1}-1$ and
  \begin{equation*}
G^{\ell-1}=
\begin{pmatrix}
G^{\ell-1}_{0} & & &\\
 & G^{\ell-1}_{1} & &\\
 & & \ddots &\\
 & & & G^{\ell-1}_{2^{\ell-1}-1}\\
\end{pmatrix},
\end{equation*}
then all the operations in this step can be represented by a large matrix-vector multiplication
 \[
 \Coef{\ell}{}=G^{\ell-1} \Coef{\ell-1}{}
 \]
  for $\ell = L/2+1, \ldots, L$.

\item{\emph{Termination.}} Finally, $\ell = L$ again and set $B=\Omega$. For each leaf interval $A=A^L_i \in T_X$, $i=0,1,\dots, 2^L-1$, use the constructed
  expansion coefficients $\delta^L$  in the last step to
  evaluate $u^B(x)$ for each $x \in A$,
  \begin{equation}
      \begin{split}
    u(x) = u^B(x)&= \sum_{t=1}^{r_\eps} \alpha^{AB}_t (x) \coef{AB}{t}\\
    &=e^{2\pi\i\Phi(x,c_B)} \sum_{t=1}^{r_\eps}\left(M_t^A(x)e^{-2\pi\i\Phi(g_t^A,c_B)} \coef{AB}{t} \right).
    \label{eqn:PIBFbf3}
      \end{split}
  \end{equation}

  For each domain pair $(A,B)=(A^L_i,\Omega)$, let $U^L_i\in\mathbb{C}^{ O(1) \times r_\eps}$ be the matrix representing the linear transformation in \eqref{eqn:PIBFbf3}. Define
  \[
  U^L=\diag\left\{U^L_{0},U^L_{1},
 \cdots,U^L_{2^{L}-1}\right\}.\]
  By the same argument in the initialization step, the matrix-vector multiplication
  format of all the operations in this step is
 \[
   u=U^L \Coef{L}{}.
   \]

\end{enumerate}

In sum, by the discussion above, $K$ is approximated by the preliminary IBF as
\[
K \approx  U^L G^{L-1}\cdots G^{L/2}  M^{L/2} (H^{L/2})^* \dots (H^{1})^* (V^0)^*.
\]
Since the total number of nonzero entries of the above sparse factors
is $O(r_\eps^2 N\log N)$ and each factors can be constructed explicitly by
the Chebyshev interpolation,
the operation and memory complexity of constructing and applying the IBF
is $O(r_\eps^2 N\log N)$.

\section{Optimal interpolative butterfly factorization}
\label{sec:optIBF}

Recall that at each level $\ell$ the kernel matrix $K$ restricted in a domain pair
$(A,B)=(A^\ell_i,B^{L-\ell}_j)$, denoted as $K^\ell_{ij}$, is a numerically low-rank matrix.
For a given $\eps$, let $r_0(\eps,\ell,i,j)$ be the numerical rank provided by a truncated
SVD of $K^\ell_{ij}$. With abuse of notation, let us use $r_0$ for simplicity. Similarly, let
$r_\eps$ denote the numerical rank provided by the low-rank approximation by
Chebyshev interpolation. Since $r_\eps$ might not be able to reveal the optimal numerical rank $r_0$, i.e. $r_\eps>r_0$. The $O(r_\eps^2)$ prefactor in the complexity
of the preliminary IBF might be far from  the optimal one, $r_0^2$.

The above observation motivates the design of the novel sweeping matrix compression based on  a sequence of structure-preserving
matrix compression. The sweeping matrix compression further
compresses the preliminary IBF into an optimal one. The main idea  is to
propagate low-rank property among matrix factors in the preliminary IBF and to
 shrink the size of dense submatrices in these matrix factors
 by the randomized low-rank approximation in Section \ref{sec:pre}.

The structure-preserving matrix compression is introduced in
Section \ref{sec:spmc}. The sweeping matrix compression
consists of two stages, the sweep-out and the sweep-in stages.
 They will be presented in Section \ref{sec:out} and Section \ref{sec:in},
respectively.

\subsection{Structure-preserving matrix compression}
\label{sec:spmc}

One key idea to construct the optimal IBF is the sweeping matrix compression via
a sequence of structure-preserving matrix compression introduced below.
Suppose $S$ is a block matrix with $m\times k$ blocks, i.e.,
\begin{equation*}
S=
  \begin{pmatrix}
      S_{0,0} & S_{0,1} & \cdots & S_{0,k-1}\\
   S_{1,0} & S_{1,1} & & S_{1,k-1}\\
      \vdots & & \ddots & \\
    S_{m-1,0} & S_{m-1,1} & &S_{m-1,k-1}
    \end{pmatrix}\in \mathbb{C}^{mr\times kr}.
\end{equation*}
For simplicity, we assume that each block in $S$ is of size $r\times r$. The structure-preserving
matrix compression can be easily extended to block matrices with different block sizes.
Let $D$ be a block-diagonal matrix with $k$ diagonal blocks and each diagonal block is of size
$r\times r_0$ with $r_0<r$, i.e.,
\begin{equation*}
D=
  \begin{pmatrix}
      D_{0} & & & \\
    & D_{1} & & \\
       & & \ddots & \\
    & & &D_{k-1}
    \end{pmatrix}\in \mathbb{C}^{kr\times kr_0}.
\end{equation*}
Let $P$ be the product of $S$ and $D$, i.e.,
\begin{equation*}
P=
\begin{pmatrix}
      P_{0,0} & P_{0,1} & \cdots & P_{0,k-1}\\
   P_{1,0} & P_{1,1} & & P_{1,k-1}\\
      \vdots & & \ddots & \\
    P_{m-1,0} & P_{m-1,1} & &P_{m-1,k-1}
    \end{pmatrix}=
  \begin{pmatrix}
      S_{0,0}D_{0}  & S_{0,1}D_{1}  & \cdots & S_{0,k-1}D_{k-1} \\
   S_{1,0}D_{0}  & S_{1,1} D_{1} & & S_{1,k-1}D_{k-1} \\
      \vdots & & \ddots & \\
    S_{m-1,0}D_{0}  & S_{m-1,1}D_{1}  & &S_{m-1,k-1}D_{k-1}
    \end{pmatrix}\in \mathbb{C}^{mr\times kr_0}.
\end{equation*}

\begin{figure}[htp]
\begin{minipage}{\textwidth}
\centering
\resizebox{4cm}{!}{
\begin{tikzpicture}[baseline=-0.5ex]
      \tikzset{every left delimiter/.style={xshift=-1ex},every right delimiter/.style={xshift=1ex}}
      \matrix (mat) [matrix of math nodes, left delimiter=(, right delimiter=)] {
      \draw;
        \draw[fill=gray] (0,8) rectangle (2,6);
        \draw[fill=gray] (4,8) rectangle (6,6);
        \draw[fill=gray] (2,6) rectangle (4,4);
         \draw[fill=gray] (4,4) rectangle (6,2);
         \draw[fill=gray] (6,4) rectangle (8,2);
         \draw[fill=gray] (2,2) rectangle (4,0);
         \draw[fill=gray] (4,2) rectangle (6,0);
         \draw (0,8) rectangle (8,0);
\\
      };
\end{tikzpicture}
}
\resizebox{1.5cm}{!}{
\begin{tikzpicture}[baseline=-0.5ex]
      \tikzset{every left delimiter/.style={xshift=-1ex},every right delimiter/.style={xshift=1ex}}
      \matrix (mat) [matrix of math nodes, left delimiter=(, right delimiter=)] {
\draw;
        \draw[fill=gray] (0,8) rectangle (0.5,6);
        \draw[fill=gray] (0.5,6) rectangle (1,4);
         \draw[fill=gray] (1,4) rectangle (1.5,2);
         \draw[fill=gray] (1.5,2) rectangle (2,0);
         \draw (0,8) rectangle (2,0);
\\
      };
\end{tikzpicture}
}
=
\resizebox{1.5cm}{!}{
\begin{tikzpicture}[baseline=-0.5ex]
      \tikzset{every left delimiter/.style={xshift=-1ex},every right delimiter/.style={xshift=1ex}}
      \matrix (mat) [matrix of math nodes, left delimiter=(, right delimiter=)] {
      \draw;
        \draw[fill=gray] (0,8) rectangle (0.5,6);
        \draw[fill=gray] (1,8) rectangle (1.5,6);
        \draw[fill=gray] (0.5,6) rectangle (1,4);
         \draw[fill=gray] (1,4) rectangle (1.5,2);
         \draw[fill=gray] (1.5,4) rectangle (2,2);
         \draw[fill=gray] (0.5,2) rectangle (1,0);
         \draw[fill=gray] (1,2) rectangle (1.5,0);
         \draw (0,8) rectangle (2,0);
\\
      };
\end{tikzpicture}
}
$\approx$
\resizebox{1.5cm}{!}{
\begin{tikzpicture}[baseline=-0.5ex]
      \tikzset{every left delimiter/.style={xshift=-1ex},every right delimiter/.style={xshift=1ex}}
      \matrix (mat) [matrix of math nodes, left delimiter=(, right delimiter=)] {
\draw;
        \draw[fill=gray] (0,8) rectangle (0.5,6);
        \draw[fill=gray] (0.5,6) rectangle (1,4);
         \draw[fill=gray] (1,4) rectangle (1.5,2);
         \draw[fill=gray] (1.5,2) rectangle (2,0);
         \draw (0,8) rectangle (2,0);
\\
      };
\end{tikzpicture}
}
\resizebox{1cm}{!}{
\begin{tikzpicture}[baseline=-0.5ex]
      \tikzset{every left delimiter/.style={xshift=-1ex},every right delimiter/.style={xshift=1ex}}
      \matrix (mat) [matrix of math nodes, left delimiter=(, right delimiter=)] {
      \draw;
        \draw[fill=gray] (0,8) rectangle (2,6);
        \draw[fill=gray] (4,8) rectangle (6,6);
        \draw[fill=gray] (2,6) rectangle (4,4);
         \draw[fill=gray] (4,4) rectangle (6,2);
         \draw[fill=gray] (6,4) rectangle (8,2);
         \draw[fill=gray] (2,2) rectangle (4,0);
         \draw[fill=gray] (4,2) rectangle (6,0);
         \draw (0,8) rectangle (8,0);
\\
      };
\end{tikzpicture}
}\\
\end{minipage}
\caption{The structure-preserving matrix compression of type $(4,4,4,1)$ $SD=P\approx \tilde{D}\tilde{S}$, where $S\in \mathbb{C}^{16\times 16}$ (left matrix), $P\in \mathbb{C}^{16\times 4}$ (middle matrix), and $\tilde{S}\in \mathbb{C}^{4\times 4}$ (right matrix) are $4\times 4$ block matrices with the same sparsity pattern, $D\in \mathbb{C}^{16\times 4}$ (middle left matrix) and $\tilde{D}\in \mathbb{C}^{16\times 4}$ (middle right matrix) are $4\times 4$ block-diagonal matrices, }
\label{fig:SPMC}
\end{figure}
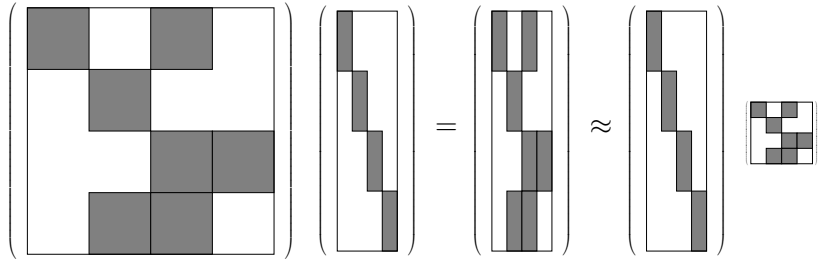

We call that the matrix pencil $(S,D)$ has the \textbf{structure-preserving
low rank property} of type $(m,k,r,r_0)$,
 if it satisfies the following condition.
 For each $i=1$, $\dots$, $m$, there exists a low-rank approximation
\begin{equation}
\label{eqn:SPMCLR}
  \begin{pmatrix}
      P_{i,0} & P_{i,1} & \cdots & P_{i,k-1}
    \end{pmatrix}
    \approx \tilde{D}_{i}
      \begin{pmatrix}
     \tilde{S}_{i,0} & \tilde{S}_{i,1} & \cdots & \tilde{S}_{i,k-1}
    \end{pmatrix}
\end{equation}
where $\tilde{D}_{i}\in \mathbb{C}^{r\times r_0}$ and
$\tilde{S}_{i,j}\in \mathbb{C}^{r_0\times r_0}$ for $j=0,1,\dots,k-1$.
Under the assumption of the structure-preserving low-rank property of type $(m,k,r,r_0)$, by assembling
the low-rank approximations in \eqref{eqn:SPMCLR} into two large factors $\tilde{D}$
and $\tilde{S}$, the structure-preserving
matrix compression of type $(m,k,r,r_0)$ factorizes $SD$ as
\[
SD\approx \tilde{D}\tilde{S},
\]
where $\tilde{D}$ is again a block-diagonal matrix with $k$ diagonal blocks  of size
$r\times r_0$, $\tilde{S}$ is a block matrix with $m\times k$ blocks of size $r_0\times r_0$,
$S$ and $\tilde{S}$ share the same column indices of nonzero blocks in each row block
(see Figure \ref{fig:SPMC} for an example).

\subsection{Sweeping matrix compression: sweep-out stage}
\label{sec:out}

The optimal IBF of $K$ is built in two stages.  In the
first stage, we further compress the matrix factors in the
preliminary IBF,
\[
K \approx  U^LG^{L-1} \cdots G^{h}   {M^{h}} (H^{h})^* \dots (H^1)^* (V^0)^*,
\]
sweeping from the middle matrix $M^h$ and moving out
towards $U^L$ and $V^0$. Recall that each nonzero submatrix
$\tilde{M}^h_{ij}\in \mathbb{C}^{r_\eps\times r_\eps}$ in $M^h$ is the kernel function $K(x.\xi)$ evaluated at the Chebyshev grid points in the domain pairs $(A,B)=(A^h_i,B^h_j)$ at level $h$. Hence,
 $\tilde{M}^h_{ij}$ is a numerically low-rank matrix if $r_\eps>r_0$ (which is usually the case in existing butterfly algorithms using Chebyshev interpolations to construct low-rank factorizations). Hence, $ {M^h}$ can be further compressed via a rank-revealing low-rank
approximation of each $\tilde{M}^h_{ij}$, e.g., the truncated SVD, into
\[
 {M^h}\approx  {C^h} {\bar{M}^h} {(R^h)}^*,
\]
resulting the {\bf middle level factorization}
\begin{equation}
\label{eqn:KeqUCMRV}
K \approx U^L G^{L-1} \cdots G^{h}   {C^h} {\bar{M}^h} {(R^h)}^* (H^{h})^* \dots (H^1)^*(V^0)^*.
\end{equation}
 The middle level
  factorization is described in detail in Section \ref{sec:MiddleConstruction}.

Next, we recursively factorize
\begin{equation*}
  {G^\ell}  {C^\ell} \approx
 {C^{\ell+1}} {\bar{G}^{\ell}}
\end{equation*}
for $\ell=h,h+1,\dots,L-1$,
\begin{equation*}
( {H^\ell}  {R^{\ell}})^* \approx ( {\bar{H}^{\ell}})^*
(  {R^{\ell-1}})^*
\end{equation*}
for $\ell=h,h-1,\dots,1$,
 since the matrix pencils $(G^\ell,C^\ell)$ and $(H^\ell, R^{\ell})$ satisfy the structure-preserving low-rank property. In another word, ${C^\ell}$ and ${R^\ell}$
propagate the low-rank property of $M^h$ to $G^\ell $ and
$H^{\ell}$, respectively, so that we can further compress $G^{\ell}$ and $H^{\ell}$. After this recursive
factorization, let
\begin{equation*}
  {\bar{U}^L}= {U^L}  {C^L},
\end{equation*}
and
\begin{equation*}
{\bar{V}^0}= {V^0}  {R^{0}},
\end{equation*}
then one reaches at a more compressed IBF of $K$:
\begin{equation}
 K \approx  \bar{U}^L\bar{G}^{L-1}\cdots \bar{G}^{h} \bar{M}^h
  (\bar{H}^{h})^* \cdots (\bar{H}^{1})^* (\bar{V}^0)^*,
\end{equation}
where all factors are sparse matrices with almost $\O(r_0^2N)$ nonzero entries.
We refer to this stage as the {\bf recursive factorization} and it is
discussed in detail in Section \ref{sec:RCU} and \ref{sec:RCV}.

\subsubsection{Middle level factorization}
\label{sec:MiddleConstruction}
Recall that we denote the $i$th node at
level $\ell$ in $T_X$ as $A^\ell_i$ for $i=0,1,\dots,2^\ell-1$ and the
$j$th node at level $L-\ell$ in $T_\Omega$ as $B^{L-\ell}_j$ for
$j=0,1,\dots,2^{L-\ell}-1$. Let $h=L/2$ and $m=2^{L/2}$.
Since the nonzero submatrix $\tilde{M}^{h}_{ij}$ in $M^h$ is a matrix representation
of the kernel function $K(x,\xi)$ evaluated at the Chebyshev grid points $\{g_t^A\}_{t}$
and  $\{g_s^B\}_{s}$ for the domain pair $(A,B)=(A^h_i,B^h_j)$ at the middle level
$\ell=h$. $\tilde{M}^{h}_{ij}\in \mathbb{C}^{r_\eps\times r_\eps}$ is numerically rank $r_0$. Hence, a
rank-$r_0$ approximation to every $\tilde{M}^{h}_{i,j}$ is computed by
the SVD algorithm via random sampling in \cite{Butterfly4,symbol3} with $O(r_\eps r_0 + r_0^3)$ operations. In fact, when $r_\eps$ is already very small, a direct method for SVD truncation of order $r_\eps^3$ is efficient as well. Once the approximate SVD
of $\tilde{M}^{h}_{i,j}$ is ready, it is transformed in the form
\[
 {\tilde{M}^{h}_{i,j} }\approx { C^h_{i,j}S^h_{i,j}(R^h_{j,i})^*}
\]
following \eqref{eq:lowrankSVD}. We would like to emphasize that the
columns of $C^h_{i,j}$ and $R^h_{j,i}$ are scaled with the singular
values of the approximate SVD so that they keep track of the
importance of these columns in approximating $\tilde{M}^{h}_{i,j}$.

After calculating the approximate rank-$r_0$ factorization of each
$\tilde{M}^{h}_{i,j}$, we assemble these factors into three block matrices
$C^h$, $\bar{M}^h$ and $R^h$ as follows:
\begin{equation}\label{eq:KeqUMV}
  \begin{split}
   {   {M}^{h} }
=&   \begin{pmatrix}
       { M^h_{0,0} }&  \cdots &  { M^h_{0,m-1}}\\
      \vdots & \ddots & \vdots \\
     {   M^h_{m-1,0} }&  \cdots  & {  M^h_{m-1,m-1}}
    \end{pmatrix}\\
    = &
    \begin{pmatrix}
     {   C^h_{0} } & &\\
      &  \ddots &\\
      & &  { C^h_{m-1}}
    \end{pmatrix}
    \begin{pmatrix}
       { \bar{M}^h_{0,0} }&  \cdots &  { \bar{M}^h_{0,m-1}}\\
      \vdots & & \ddots &  \vdots \\
     {   \bar{M}^h_{m-1,0} }& \cdots & { \bar{M}^h_{m-1,m-1}}
    \end{pmatrix}
    \begin{pmatrix}
   {     (R_{0}^h)^* }&  &\\
      &  \ddots &\\
      &  & {  (R_{m-1}^h)^*}
    \end{pmatrix}\\
    = &  {C^h \bar{M}^h (R^h)^*  },\\
  \end{split}
\end{equation}
where
\begin{equation}
 {  C_{i}^h}=
     \begin{pmatrix}
    {   C^h_{i,0} }& & &\\
      & { C^h_{i,1} }& &\\
      & & \ddots &\\
      & & & { C^h_{i,m-1}}
    \end{pmatrix}
  \in \bbC^{mr_\eps \times mr_0},
  \label{eq:expression-C}
  \end{equation}
  \begin{equation}
 {  R^h_{j}}=
 \begin{pmatrix}
    {   R^h_{0,j} }& & &\\
      & { R^h_{1,j} }& &\\
      & & \ddots &\\
      & & & { R^h_{m-1,j}}
    \end{pmatrix}
  \in \bbC^{mr_0 \times mr_\eps},
  \label{eq:expression-R}
\end{equation}
and $\bar{M}^{h}_{i,j}\in \mathbb{C}^{mr_0\times mr_0}$ is also a $m\times m$ block matrix with block size $r_0\times r_0$ where all blocks are zero except that the $(j,i)$ block is equal to the diagonal matrix ${S}^h_{i,j}\in \mathbb{C}^{r_0\times r_0}$.

\begin{figure}[htp]
\begin{minipage}{\textwidth}
\centering
\resizebox{5.5cm}{!}{
\begin{tikzpicture}[baseline=-0.5ex]
      \tikzset{every left delimiter/.style={xshift=-1ex},every right delimiter/.style={xshift=1ex}}
      \matrix (mat) [matrix of math nodes, left delimiter=(, right delimiter=)] {
      \draw;
\foreach \i in {0,1,...,3}{
    \foreach \j in {0,1,...,3}{
        \draw[fill=gray] (4*\i+\j,8-4*\j-\i) rectangle (4*\i+\j+1,7-4*\j-\i);
    }
}
         \draw (0,8) rectangle (16,-8);
\\
      };
\end{tikzpicture}
}
$\approx$
\resizebox{1.5cm}{!}{
\begin{tikzpicture}[baseline=-0.5ex]
      \tikzset{every left delimiter/.style={xshift=-1ex},every right delimiter/.style={xshift=1ex}}
      \matrix (mat) [matrix of math nodes, left delimiter=(, right delimiter=)] {
\draw;
\foreach \i in {0,1,...,15}{
        \draw[fill=gray] (0+\i*0.5,32-2*\i) rectangle (0.5+\i*0.5,30-2*\i);
         }
         \draw (0,32) rectangle (8,0);
\\
      };
\end{tikzpicture}
}
\resizebox{1.5cm}{!}{
\begin{tikzpicture}[baseline=-0.5ex]
      \tikzset{every left delimiter/.style={xshift=-1ex},every right delimiter/.style={xshift=1ex}}
      \matrix (mat) [matrix of math nodes, left delimiter=(, right delimiter=)] {
      \draw;
\foreach \i in {0,1,...,3}{
    \foreach \j in {0,1,...,3}{
        \draw[fill=gray] (4*\i+\j,8-4*\j-\i) rectangle (4*\i+\j+1,7-4*\j-\i);
    }
}
         \draw (0,8) rectangle (16,-8);
\\
      };
\end{tikzpicture}
}
\resizebox{5.5cm}{!}{
\begin{tikzpicture}[baseline=-0.5ex]
      \tikzset{every left delimiter/.style={xshift=-1ex},every right delimiter/.style={xshift=1ex}}
      \matrix (mat) [matrix of math nodes, left delimiter=(, right delimiter=)] {
\draw;
\foreach \i in {0,1,...,15}{
        \draw[fill=gray] (0+\i*2,8-0.5*\i) rectangle (2+\i*2,7.5-0.5*\i);
         }
         \draw (0,8) rectangle (32,0);
\\
      };
\end{tikzpicture}
}\\
\end{minipage}
\caption{The middle level factorization of
    a $64\times 64$ matrix ${M}^2\approx C^2\bar{M}^2(R^2)^*$
    assuming $r_0=1$ and $r_\eps=4$.
    Grey blocks indicate nonzero blocks.
    $C^2$ and $R^2$ are block-diagonal matrices with $16$ blocks of size $4 \times 1$.
    The diagonal blocks of $C^2$ and $R^2$ are assembled
    according to Equation \eqref{eq:expression-C} and \eqref{eq:expression-R}, respectively, 
    as indicated by black rectangles.
    $\bar{M}^2$ is a $4\times 4$ block matrix with each block
    $\bar{M}^2_{i,j}$ itself an $4\times 4$ block matrix
    containing diagonal weights matrix on the $(j,i)$ block.}
\label{fig:compression-1D}
\end{figure}
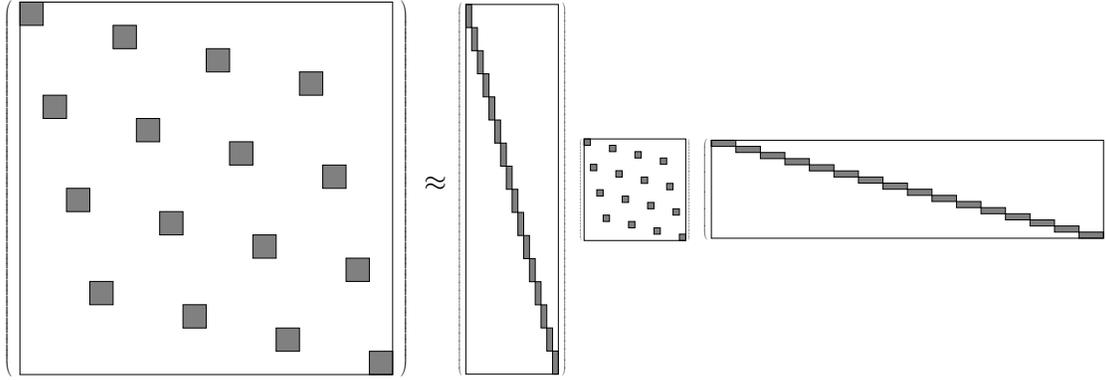

It is obvious that there are only $r_0N$ nonzero entries in $\bar{M}^h$ and $r_\eps r_0N$ nonzero entries in $C^h$ and $R^h$.
See Figure~\ref{fig:compression-1D} for an example of a middle level
factorization of a $64\times 64$ matrix with $r_0=1$ and $r_\eps=4$.

%

\subsubsection{Recursive factorization towards $U^L$}
\label{sec:RCU}

Each recursive factorization at level $\ell$
\begin{equation}
 \label{eqn:GC}
  {G^\ell}  {C^\ell} \approx
 {C^{\ell+1}} {\bar{G}^{\ell}}
\end{equation}
 results from the structure-preserving low-rank property that originates from the
low-rank property of $K^\ell_{i,j}$ for all index pairs $(i,j)$.  At level
$\ell$, recall that
  \begin{equation*}
G^{\ell}=
\begin{pmatrix}
G^{\ell}_{0} & & &\\
 & G^{\ell}_{1} & &\\
 & & \ddots &\\
 & & & G^{\ell}_{2^{\ell}-1}\\
\end{pmatrix},
\end{equation*}
where
 \begin{equation*}
     G^{\ell}_{i}=
     \begin{pmatrix}
     G^{\ell}_{2i,0} & G^{\ell}_{2i,1} & & &\\
      & & G^{\ell}_{2i,2} & G^{\ell}_{2i,3} & &\\
      & & & & \ddots &\\
      & & & & & G^{\ell}_{2i,2^{L-\ell-1}-2} & G^{\ell}_{2i,2^{L-\ell-1}-1}\\
     \midrule
     G^{\ell}_{2i+1,0} & G^{\ell}_{2i+1,1} & & &\\
      & & G^{\ell}_{2i+1,2} & G^{\ell}_{2i+1,3} & &\\
      & & & & \ddots &\\
      & & & & & G^{\ell}_{2i+1,2^{L-\ell-1}-2}
       & G^{\ell}_{2i+1,2^{L-\ell-1}-1}
     \end{pmatrix}
\end{equation*}
for $i=0,1,\dots,2^{\ell}-1$.
By construction, $[ G^{\ell}_{2i,2j}\quad G^{\ell}_{2i,2j+1}]$ is an
interpolation matrix that interpolates the kernel function $K(x,\xi)$
from Chebyshev grid points in the domain pairs $(A^{\ell}_i,B^{L-\ell}_{2j})$
and $(A^{\ell}_i,B^{L-\ell}_{2j+1})$ to the points
in the domain pair $(A^{\ell+1}_{2i},B^{L-\ell-1}_{j})$.
Similarly, $[G^{\ell}_{2i+1,2j}\quad G^{\ell}_{2i+1,2j+1}]$ interpolates the kernel function $K(x,\xi)$
from $(A^{\ell}_i,B^{L-\ell}_{2j})$
and $(A^{\ell}_i,B^{L-\ell}_{2j+1})$
to $(A^{\ell+1}_{2i+1},B^{L-\ell-1}_{j})$.

At level $\ell=h$,
\begin{equation*}
C^h=
    \begin{pmatrix}
    {   C^h_{0} } & &\\
      &  \ddots &\\
      & & { C^h_{m-1}}
    \end{pmatrix}
\end{equation*}
where
\begin{equation*}
 {  C_{i}^h}=
     \begin{pmatrix}
    {   C^h_{i,0} }& & &\\
      & { C^h_{i,1} }& &\\
      & & \ddots &\\
      & & & { C^h_{i,m-1}}
    \end{pmatrix}
  \in \bbC^{r_\eps m \times r_0 m}.
  \end{equation*}
By construction, the column space of $C^h_{i,j}$ comes from the column space of $K^h_{i,j}$ for all $i=0,1,\dots,m-1$ and $j=0,1,\dots,m-1$. Hence,
$
\begin{pmatrix}
G^{h}_{2i,2j}C^h_{i,2j} & G^{h}_{2i,2j+1}C^h_{i,2j+1}
\end{pmatrix}
$
represents the column space of $K^{h+1}_{2i,j}$, which implies that
$
\begin{pmatrix}
G^{h}_{2i,2j}C^h_{i,2j} & G^{h}_{2i,2j+1}C^h_{i,2j+1}
\end{pmatrix}
$
is numerically rank $r_0$. By the randomized low-rank approximation, we have
\[
\begin{pmatrix}
G^{h}_{2i,2j}C^h_{i,2j} & G^{h}_{2i,2j+1}C^h_{i,2j+1}
\end{pmatrix}
\approx
C^{h+1}_{2i,j}
\begin{pmatrix}
\bar{G}^{h}_{2i,2j} & \bar{G}^{h}_{2i,2j+1}
\end{pmatrix},
\]
where $C^{h+1}_{2i,j} \in \mathbb{C}^{r_\eps\times r_0}$,
$\bar{G}^{h}_{2i,2j} \in \mathbb{C}^{r_0\times r_0}$,
and $\bar{G}^{h}_{2i,2j+1} \in \mathbb{C}^{r_0\times r_0}$.
By a similar argument, we have
\[
\begin{pmatrix}
G^{h}_{2i+1,2j}C^h_{i,2j} & G^{h}_{2i+1,2j+1}C^h_{i,2j+1}
\end{pmatrix}
\approx
C^{h+1}_{2i+1,j}
\begin{pmatrix}
\bar{G}^{h}_{2i+1,2j} & \bar{G}^{h}_{2i+1,2j+1}
\end{pmatrix},
\]
where $C^{h+1}_{2i+1,j} \in \mathbb{C}^{r_\eps\times r_0}$,
$\bar{G}^{h}_{2i+1,2j} \in \mathbb{C}^{r_0\times r_0}$,
and $\bar{G}^{h}_{2i+1,2j+1} \in \mathbb{C}^{r_0\times r_0}$.
Hence, the matrix pencil $(G^h,C^h)$ satisfies
the conditions of the structure-preserving low-rank property of type
$(2^L,2^L,r_\eps,r_0)$.
Assembling $\{\bar{G}^{h}_{i,j}\}$ and $\{C^{h+1}_{i,j}\}$
together to generate $\bar{G}^h\in \mathbb{C}^{r_0N\times r_0N}$
and $C^{h+1}\in \mathbb{C}^{r_\eps N\times r_0 N}$, respectively,
in the same way as generating $G^h$ and $C^h$, we have
\[
  {G^h}  {C^h} \approx
 {C^{h+1}} {\bar{G}^{h}}.
\]

At other levels $\ell=h,\dots,L-1$, similarly to the discussion above,
the matrix pencil $(G^\ell,C^\ell)$ satisfies
the structure-preserving low-rank property of type $(2^L,2^L,r_\eps,r_0)$.
By assembling the results of randomized low-rank approximations, we have
\[
  {G^\ell}  {C^\ell} \approx
 {C^{\ell+1}} {\bar{G}^{\ell}}
\]
for $\ell=h,\dots,L-1$.

At the final level $\ell=L$, both matrices $U^L,C^L$
are block diagonal matrices, we simply multiply them together and let
\[
   \bar{U}^L={U^L}  {C^L} .
\]

After the recursive factorization of each $G^\ell C^\ell$, we have
\begin{equation}
 K \approx  \bar{U}^L\bar{G}^{L-1}\cdots \bar{G}^{h} \bar{M}^h
 ({R}^{h})^* ({H}^{h})^* \cdots({H}^{1})^* ({V}^{0})^*,
\end{equation}
where $\bar{G}^\ell$ contains only $r_0^2 N$ nonzero entries and $\bar{U}^{L}$ contains $r_0 N$ nonzero entries.

\subsubsection{Recursive factorization towards $V^0$}
\label{sec:RCV}

The recursive factorization towards $V^0$ is similar to the one towards $U^L$.
In each step of the factorization
\begin{equation}
 \label{eqn:HR}
( {H^\ell}  {R^{\ell}})^* \approx ( {\bar{H}^{\ell}})^*
(  {R^{\ell-1}})^*
\end{equation}
for all $\ell=h,h+1,\dots,1$,
we take advantage of the structure-preserving low-rank property of
the matrix pencil $(H^\ell,R^\ell)$.
The same procedure of Section
\ref{sec:RCU} now to $(H^\ell,R^\ell)$ leads to the recursive factorization
 in \eqref{eqn:HR}. Define
\begin{equation}
\label{eqn:V0}
\bar{V}^{0} = {V^0}  {R^{0}},
\end{equation}
then we have
\begin{equation*}
 ({R}^{h})^* ({H}^{h})^* \cdots({H}^{1})^* ({V}^{0})^*\approx
 { (\bar{H}^{h})^* \cdots (\bar{H}^{1})^* (\bar{V}^{0})^* },
\end{equation*}
where $\bar{H}^\ell$ contains only $r_0^2 N$ nonzero entries and $\bar{V}^0$ contains $r_0 N$ nonzero entries.

After the recursive factorization sweeping from the middle matrix towards $U^L$ and $V^0$,  we reach a more compressed IBF
\begin{equation}
\label{eqn:IBFout}
 K \approx  \bar{U}^{L} \bar{G}^{L-1}\cdots \bar{G}^{h} \bar{M}^h
 (\bar{H}^{h})^* \cdots (\bar{H}^{1})^*  (\bar{V}^{0})^*.
\end{equation}

\subsection{Sweeping matrix compression: sweep-in stage}
\label{sec:in}

If the points in the sets $X$ and $\Omega$ are distributed uniformly,
the IBF in \eqref{eqn:IBFout} is already optimal
in the butterfly factorization scheme,
i.e., nearly all dense submatrices in its factors are
of size $r_0\times r_0$, where $r_0$ is the numerical rank of
the kernel function $K(x,\xi)$ sampled uniformly
in a domain pair $(A,B)\in T_X\times T_\Omega$.
However, when the point sets are nonuniform,
e.g., in the nonuniform FFT,
the number of samples in $(A,B)$ might be far smaller than $r_0$.
This means that there might be dense submatrices of size
less than $r_0\times r_0$ in the block diagonal matrices $ \bar{U}^{L}$
and $ \bar{V}^{0}$.
This motivates a sequence of structure-preserving low-rank
matrix compression to further compress the IBF in \eqref{eqn:IBFout},
sweeping from outer matrices and moving
towards the middle matrix $\bar{M}^h$ as follows.
\begin{equation*}
\bar{U}^{L} \approx \dot{U}^{L} C^{L}
\end{equation*}
and
\begin{equation*}
C^{\ell+1}  \bar{G}^\ell \approx \dot{G}^\ell C^{\ell}
\end{equation*}
for $\ell=L-1,L-2,\dots,h$;
similarly, we have,
\begin{equation*}
(\bar{V}^{0})^* \approx (R^{0})^* (\dot{V}^{0})^*
\end{equation*}
and
\begin{equation*}
(\bar{H}^\ell)^* (R^{\ell-1})^* \approx (R^{\ell})^* (\dot{H}^\ell)^*
\end{equation*}
for $\ell=1,2,\dots,h$.
The sweeping matrix compression above is due to the fact that
the matrix pencils $((\bar{G}^\ell)^*, (C^{\ell+1})^*)$
and $( (\bar{H}^\ell)^*,(R^{\ell-1})^*)$ satisfy
the structure-preserving low-rank property.
In another word, ${C^\ell}$ and ${R^\ell}$
propagate the low-rank property of $\bar{U}^L$ and $\bar{V}^L$ to
$\bar{G}^\ell$ and
$\bar{H}^{\ell}$, respectively,
so that we can further compress $\bar{G}^{\ell}$ and $\bar{H}^{\ell}$.
After this recursive
factorization, let $\dot{M}^h = C^h \bar{M}^h(R^h)^*$,
then one reaches the optimal IBF of $K$:
\begin{equation}
\label{eqn:IBFin}
K \approx \dot{U}^{L}\dot{G}^{L-1}\cdots\dot{G}^{h}\dot{M}^h
(\dot{H}^{h})^* \cdots (\dot{H}^{1})^* (\dot{V}^{0})^*.
\end{equation}
where all factors are sparse matrices with almost $\O(r_0^2N)$
or less nonzero entries.

\vspace{10 mm}
For a given input vector $g\in\bbC^N$, the $\O(N^2)$ matrix-vector
multiplication $u = Kg$ can be approximated by a sequence of $\O( \log
N)$ sparse matrix-vector multiplications given by the optimal IBF.
Since computing the factors in the preliminary IBF takes
only $O(r_\eps^2 N\log N)$ operations
and there are only $O(r_0^2 N\log N)$ nonzero entries in the optimal IBF,
the construction and application complexity in operation is
$O(r_\eps^2 N\log N)$ and $O(r_0^2 N\log N)$, respectively.
Since all the matrix factors in the preliminary IBF
can be generated explicitly,
the peak memory complexity $O(r_\eps^2 N\log N)$ occurs
when the preliminary IBF is completed.

\section{High dimensional extension}
\label{sec:mIBF}

Although we limited our discussion to one-dimensional problems in previous discussion,
the interpolative butterfly factorization, along with its construction algorithm, can be
easily generalized to higher dimensions. 

By the theorems in \cite{FIO09,FIO14}, a multidimensional kernel function
$K(x,\xi)$ satisfying Assumption \ref{asm} is
complementary low-rank (e.g., the nonuniform
FFT). In this case, similarly to Section \ref{sec:preIBF},
by writing the multidimensional butterfly algorithm
\cite{FIO09,FIO14} into a matrix factorization form,
we have a preliminary multidimensional
IBF
\[
K \approx  U^L G^{L-1}\cdots G^{h}  M^{h} (H^{h})^* \dots (H^{1})^* (V^0)^*,
\]
where the depth $L=\O(\log N)$ of $T_X$ and $T_\Omega$ is assumed to be even,
$h=L/2$ is a middle level index, all factors are sparse matrices with
$\O(N)$ nonzero entries and a large prefactor.
The preliminary IBF can be further compressed by the sweeping
matrix compression to obtain the optimal IBF
\[
K \approx \dot{U}^{L} \dot{G}^{L-1} \cdots \dot{G}^{h} \dot{M}^h
(\dot{H}^{h})^* \cdots (\dot{H}^{1})^* (\dot{V}^{0})^*.
\]

However, many important
multidimensional kernel matrices fail to satisfy Assumption \ref{asm}
and is not complementary low-rank in the entire domain $X \times \Omega$.
The most significant example,
the multidimensional Fourier integral operator,
typically has a singularity when $\xi=0$ in the
$\Omega$ domain. Fortunately, it was proved that this kind of kernel functions
satisfies the complementary low-rank property in the domain away from $\xi=0$.
An {\em multiscale interpolative butterfly factorization} (MIBF)
hierarchically partitions the domain $\Omega$ into subdomains $\{\Omega_t\}_t$ excluding
the singular point $\xi=0$ and apply the multidimensional IBF to the kernel restricted in
each subdomain pair $X\times \Omega_t$.

To be more specific, in two-dimensional problems, suppose
\begin{equation*}
  \Omega = \left\{ \xi =(n_1,n_2),- \frac{N^{1/2}}{2} \leq n_1,n_2 < \frac{N^{1/2}}{2}\text{ with }
  n_1,n_2 \in \mathbb{Z} \right\}.
\end{equation*}
Let
\begin{equation}\label{eq:domaindecomp}
  \Omega_t =
  \left\{(\xi_1,\xi_2):\frac{N^{1/2}}{2^{t+2}}<\max(|\xi_1|,|\xi_2|) \leq
  \frac{N^{1/2}}{2^{t+1}}\right\}\cap \Omega,
\end{equation}
for $t=0,1,\dots,\log_2 n-s$, where $n=N^{1/2}$ and $s$ is a small constant,
and $\Omega_C=\Omega\setminus\cup_t\Omega_t$. Equation
\eqref{eq:domaindecomp} is a corona decomposition of $\Omega$, where
each $\Omega_t$ is a corona subdomain and $\Omega_C$ is a square
subdomain at the center containing $\O(1)$ points.

The FIO kernel function satisfies the complementary
low-rank property when it is restricted in each subdomain $X\times
\Omega_t$ as proved in \cite{FIO14}. Hence, the
MIBF evaluates
\[
u(x) = \sum_{\xi\in\Omega} e^{2\pi\imath \Phi(x,\xi)}g(\xi)
\]
via a multiscale summation,
\begin{equation}\label{eq:mba}
  u(x) = u_C(x) + \sum_{t=0}^{\log_2n-s} u_t(x) = \sum_{\xi\in
    \Omega_C}e^{2\pi\imath \Phi(x,\xi)}g(\xi)
  +\sum_{t=0}^{\log_2n-s}\sum_{\xi\in \Omega_t}e^{2\pi\imath
    \Phi(x,\xi)}g(\xi).
\end{equation}
For each $t$, the MIBF evaluates $u_t(x) = \sum_{\xi\in\Omega_t}
e^{2\pi\imath\Phi(x,\xi)} g(\xi)$ with a standard multidimensional IBF
and the final piece
$u_C(x)$ is evaluated directly in $\O(N)$ operations.
Let $K_t$ and $K_C$ denote the matrix representation of the kernel $K(x,\xi)$ restricted
in $X\times \Omega_t$ and $X\times \Omega_C$, respectively,
then by the standard multidimensional IBF, we have
  \begin{equation*}
    K_t \approx \dot{U}_t^{L_t}\dot{G}_t^{L_t-1}\cdots
    \dot{G}_t^{\frac{L_t}{2}}
    \dot{M}_t^{\frac{L_t}{2}}
    \left(\dot{H}_t^{\frac{L_t}{2}}\right)^*\cdots
    \left(\dot{H}_t^{1}\right)^*\left(\dot{V}_t^{0}\right)^*.
  \end{equation*}
  Once we have computed the optimal IBF in each restricted domain,
  the multiscale summation in \eqref{eq:mba} is approximated by
  \begin{equation}
   u= Kg \approx K_C R_Cg + \sum_{t=0}^{\log_2 n-s}
    \dot{U}_t^{L_t} \dot{G}_t^{L_t-1} \cdots \dot{M}_t^{\frac{L_t}{2}}
    \cdots
    \left(\dot{H}_t^{1}\right)^*\left(\dot{V}_t^{0}\right)^* R_t g,
    \label{eq:f1}
  \end{equation}
  where $R_C$ and $R_t$ are the
  restriction operators to the domains $\Omega_C$ and $\Omega_t$
  respectively.

 The construction and application complexity of the MIBF
 is $O(N\log N)$ with an optimally small prefactor in the butterfly scheme.

\section{Numerical results}
\label{sec:results}

This section presents several numerical examples to demonstrate the efficiency of
the interpolative butterfly factorization.
The numerical results were obtained
on a single node of a server cluster with
quad socket Intel\textsuperscript{\textregistered}
Xeon\textsuperscript{\textregistered}
CPU E5-4640 @ 2.40GHz (8 core/socket)
and 1.5 TB RAM.
All implementations are in MATLAB\textsuperscript{\textregistered}
and can found in authors' homepages.

Let $\{u^d(x),x\in X\}$, $\{u^i(x),x\in X\}$
and $\{u^m(x),x\in X\}$ denote the results
given by the direct matrix-vector multiplication,
the interpolative butterfly factorization (IBF)
and the multiscale interpolative butterfly factorization (MIBF).
The accuracies of applying the IBF and MIBF are estimated
by the relative error defined as follows,
\begin{equation}
\epsilon^i = \sqrt{\cfrac{\sum_{x\in S}|u^i(x)-u^d(x)|^2}
{\sum_{x\in S}|u^d(x)|^2}}
\quad \text{and} \quad
\epsilon^m = \sqrt{\cfrac{\sum_{x\in S}|u^m(x)-u^d(x)|^2}
{\sum_{x\in S}|u^d(x)|^2}},
\label{eq:relativeerr}
\end{equation}
where $S$ is a point set of size $256$ randomly sampled from $X$.
Meanwhile, let $R_{comp}$ denotes the compression ratio of the optimal
interpolative butterfly factorization
against the preliminary interpolative butterfly factorization,
which is defined as
\begin{equation}
R_{comp} =
\frac{\text{Memory usage of the preliminary
interpolative butterfly factorization}}
{\text{Memory usage of the optimal interpolative butterfly factorization}}.
\label{eq:compratio}
\end{equation}
$R_{comp}$ accurately evaluates the rank compression ratio
in the optimal interpolative butterfly factorization
without lost of accuracy.

Although we define $R_{comp}$ as a ratio of memory usage, it also reflects
the ratio of running time. Since the application of IBF is a sequence of
matrix-vector multiplications, the total running time is linearly depends
on the number of non-zeros which is linearly depends on the memory usage.
Therefore, $R_{comp}$ also equals the running time of the preliminary IBF
over the optimal IBF.

\subsection{Nonuniform Fourier Transform in 1D}
We first present the numerical result of the most widely used
Fourier integral operator, Fourier transform.
More specifically, we focus on one-dimensional
nonuniform Fourier transform of type I
\begin{equation}\label{eqn:1D-FT}
\widehat{u}(\xi_i) = \sum_{x_j}e^{-2\pi \imath x_j\xi_i}u(x_j),
\quad i,j=1,2,\dots,N,
\end{equation}
where $\{x_i\}$ are $N$ random points in $[0,1)$
each of which is drawn from uniform distribution $[0,1)$, and
$\xi_j = j-1-N/2$.
The values of the input function $\{{u}(x_j)\}_{j=1}^N$,
are randomly generated.

Table~\ref{tab:1D-FT} summarizes the results of this example for varying
problem sizes, $N$, and numbers of Chebyshev points, $r_\epsilon$.  In
Table~\ref{tab:1D-NUFFT}, we further provide the detailed comparison between
the application cost of IBF and unifom/non-uniform FFTs.

\begin{table}[htp]
\centering
\begin{tabular}{rcccccccc}
\toprule
   $N,r_\epsilon$ & $\epsilon^i$ & $R_{comp}$ & $T_{Factor}(min)$
                  & $T_d(sec)$ & $T_{app}(sec)$ & $T_d/T_{app}$\\
\toprule
    256, 6 & 4.35e-04 & 1.33 & 1.49e-02 & 2.73e-02 & 6.99e-04 & 3.90e+01 \\
   1024, 6 & 7.80e-04 & 1.38 & 9.47e-02 & 1.82e-01 & 1.78e-03 & 1.03e+02 \\
   4096, 6 & 8.89e-04 & 1.40 & 5.20e-01 & 1.48e+00 & 7.58e-03 & 1.96e+02 \\
  16384, 6 & 1.09e-03 & 1.42 & 2.66e+00 & 1.25e+01 & 3.66e-02 & 3.41e+02 \\
  65536, 6 & 1.12e-03 & 1.42 & 1.31e+01 & 1.60e+02 & 1.68e-01 & 9.54e+02 \\
 262144, 6 & 1.20e-03 & 1.43 & 6.16e+01 & 2.69e+03 & 7.63e-01 & 3.53e+03 \\
1048576, 6 & 1.18e-03 & 1.43 & 2.66e+02 & 4.30e+04 & 3.56e+00 & 1.21e+04 \\
\toprule
    256,10 & 3.57e-08 & 1.50 & 1.56e-02 & 2.66e-02 & 8.67e-04 & 3.07e+01 \\
   1024,10 & 5.09e-08 & 1.44 & 1.04e-01 & 1.85e-01 & 3.58e-03 & 5.17e+01 \\
   4096,10 & 1.02e-07 & 1.46 & 5.76e-01 & 1.54e+00 & 1.55e-02 & 9.92e+01 \\
  16384,10 & 1.13e-07 & 1.49 & 2.95e+00 & 1.27e+01 & 7.59e-02 & 1.67e+02 \\
  65536,10 & 1.27e-07 & 1.53 & 1.45e+01 & 1.75e+02 & 3.36e-01 & 5.22e+02 \\
 262144,10 & 1.34e-07 & 1.55 & 6.86e+01 & 2.57e+03 & 1.57e+00 & 1.64e+03 \\
1048576,10 & 1.43e-07 & 1.56 & 2.99e+02 & 4.26e+04 & 1.44e+01 & 2.96e+03 \\
\bottomrule
\end{tabular}
\caption{Numerical results for the one-dimensional Fourier transform given
  in \eqref{eqn:1D-FT}. $N$ is the problem size; $r_\epsilon$ is the number
  of Chebyshev points; $\epsilon^i$ is the sampled relative error given in
  \eqref{eq:relativeerr}; $R_{comp}$ is the compression ratio defined as
  \eqref{eq:compratio}; $T_{Factor}$ is the construction time of the IBF;
  $T_{d}$ is the running time of the direct evaluation; $T_{app}$ is the
  application time of the IBF; $T_d/T_{app}$ is the speedup factor over the
  direct evaluation. }
\label{tab:1D-FT}
\end{table}

\begin{table}[htp]
\centering
\begin{tabular}{rcccccccc}
\toprule
    $N$ & $\epsilon^i$ & $T_{app}(sec)$ & $P_{op}$ & $T_{NUFFT}(sec)$ &
    $T_{app}/T_{NUFFT}$ \\
\toprule
    256 & 4.35e-04 & 6.99e-04 & 6.64e+00 & 2.25e-04 & 3.11e+00 \\
   1024 & 7.80e-04 & 1.78e-03 & 7.82e+00 & 6.17e-04 & 2.88e+00 \\
   4096 & 8.89e-04 & 7.58e-03 & 8.65e+00 & 2.57e-03 & 2.95e+00 \\
  16384 & 1.09e-03 & 3.66e-02 & 9.24e+00 & 9.88e-03 & 3.70e+00 \\
  65536 & 1.12e-03 & 1.68e-01 & 9.71e+00 & 4.13e-02 & 4.07e+00 \\
 262144 & 1.20e-03 & 7.63e-01 & 1.01e+01 & 1.80e-01 & 4.25e+00 \\
1048576 & 1.18e-03 & 3.56e+00 & 1.04e+01 & 7.74e-01 & 4.60e+00 \\
\toprule
    256 & 3.57e-08 & 8.67e-04 & 1.41e+01 & 2.18e-04 & 3.97e+00 \\
   1024 & 5.09e-08 & 3.58e-03 & 1.93e+01 & 6.34e-04 & 5.65e+00 \\
   4096 & 1.02e-07 & 1.55e-02 & 2.23e+01 & 2.59e-03 & 5.97e+00 \\
  16384 & 1.13e-07 & 7.59e-02 & 2.44e+01 & 1.00e-02 & 7.57e+00 \\
  65536 & 1.27e-07 & 3.36e-01 & 2.57e+01 & 4.11e-02 & 8.18e+00 \\
 262144 & 1.34e-07 & 1.57e+00 & 2.66e+01 & 1.74e-01 & 9.01e+00 \\
1048576 & 1.43e-07 & 1.44e+01 & 2.74e+01 & 7.45e-01 & 1.93e+01 \\
\bottomrule
\end{tabular}
\caption{Numerical comparison between IBF and NUFFT~\cite{NUFFT} for the
one-dimensional Fourier transform given in \eqref{eqn:1D-FT}.  $P_{op}$ is
the operator-wise penalty over the uniform FFT \cite{johnson2007modified,
bernstein2007tangent}, i.e., the number of operation count over the one of
the FFT; $T_{NUFFT}$ is the running time of the NUFFT~\cite{NUFFT} where the
implementation is in Fortran; $T_{app}/T_{NUFFT}$ is the time penalty of the
non-uniform FFT. }
\label{tab:1D-NUFFT}
\end{table}

The result in Table~\ref{tab:1D-FT} reflects the $\O(N\log N)$ complexity
for both the construction and application. The relative error slightly
increases as the problem size $N$ increases. In general, nonuniform FFT
requires more effect to interpolate the irregular point distribution, which
means there is an underlying penalty factor coming from the problem itself.
Based on the numbers in Table~\ref{tab:1D-FT}, the penalty factor is on
average 9 for approximation with accuracy 1e-3 and 25 for approximation
with accuracy 1e-7. This implies that if the proposed algorithm is
well implemented and the code is deeply optimized, the application time
of the IBF for the nonuniform Fourier transform is about $9$ and $25$
times slower than the FFT for an approximation accuracy 1e-3 and 1e-7,
respectively.  Table~\ref{tab:1D-NUFFT} provides the concrete running
time comparison. The actual time penalty over the NUFFT~\cite{NUFFT}
is on average about $3$ and $6$. As we shall discuss later, the IBF would have 
better scalability in distributed and parallel computing (future work) than the 
existing NUFFT framework. Hence, the distributed and parallel IBF 
could be better than the existing NUFFT framework for large-scale computing.

\subsection{General Fourier integral operator in 1D}
Our second example is to evaluate a one-dimensional FIO \cite{BF} of the
following form:
\begin{equation}\label{eq:example1}
u(x) = \int_{\mathbb{R}}e^{2\pi \imath \Phi(x,\xi)}\widehat{f}(\xi)d\xi,
\end{equation}
where $\widehat{f}$ is the {Fourier} transform of $f$,
and $\Phi(x,\xi)$ is a phase function given by
\begin{equation}
\Phi(x,\xi) = x\cdot \xi + c(x)|\xi|,~~~c(x) = (2+\sin(2\pi x))/8.
\label{eqn:1D-FIO-kernal}
\end{equation}
The discretization of \eqref{eq:example1} is
\begin{equation}
\label{eqn:1D-FIO}
u(x_i) = \sum_{\xi_j}e^{2\pi \imath \Phi(x_i,\xi_j)}\widehat{f}(\xi_j),
\quad i,j=1,2,\dots,N,
\end{equation}
where $\{x_i\}$ and $\{\xi_j\}$ are points uniformly distributed
in $[0,1)$ and $[-N/2,N/2)$ following
\begin{equation}\label{eqn:1D-xandxi}
x_i = (i-1)/N \text{ and } \xi_j = j-1-N/2.
\end{equation}

Table~\ref{tab:1D-FIO} summarizes the results of
this example for different grid sizes $N$ and Chebyshev points $r_\epsilon$.

\begin{table}[htp]
\centering
\begin{tabular}{rcccccc}
\toprule
   $N,r_\epsilon$ & $\epsilon^i$ & $R_{comp}$ & $T_{Factor}(min)$
                             & $T_d(sec)$ & $T_{app}(sec)$ & $T_d/T_{app}$\\
\toprule
    256, 7 & 4.58e-03 & 2.19 & 1.52e-02 & 2.86e-02 & 8.26e-04 & 3.47e+01 \\
   1024, 7 & 6.53e-03 & 2.28 & 9.38e-02 & 1.84e-01 & 1.78e-03 & 1.03e+02 \\
   4096, 7 & 7.68e-03 & 2.34 & 5.05e-01 & 1.47e+00 & 8.57e-03 & 1.71e+02 \\
  16384, 7 & 8.22e-03 & 2.38 & 2.57e+00 & 1.23e+01 & 2.82e-02 & 4.37e+02 \\
  65536, 7 & 1.04e-02 & 2.41 & 1.25e+01 & 1.48e+02 & 1.23e-01 & 1.20e+03 \\
 262144, 7 & 1.05e-02 & 2.45 & 5.91e+01 & 2.48e+03 & 5.93e-01 & 4.18e+03 \\
1048576, 7 & 1.25e-02 & 2.50 & 2.59e+02 & 5.70e+04 & 2.39e+00 & 2.39e+04 \\
\toprule
    256,10 & 1.87e-05 & 1.82 & 1.60e-02 & 2.78e-02 & 9.81e-04 & 2.84e+01 \\
   1024,10 & 9.47e-06 & 1.87 & 9.99e-02 & 1.86e-01 & 3.08e-03 & 6.03e+01 \\
   4096,10 & 1.03e-05 & 2.00 & 5.48e-01 & 1.50e+00 & 1.19e-02 & 1.26e+02 \\
  16384,10 & 1.09e-05 & 2.07 & 2.80e+00 & 1.22e+01 & 5.76e-02 & 2.12e+02 \\
  65536,10 & 1.29e-05 & 2.14 & 1.37e+01 & 1.51e+02 & 3.09e-01 & 4.88e+02 \\
 262144,10 & 1.37e-05 & 2.18 & 6.45e+01 & 2.58e+03 & 1.13e+00 & 2.28e+03 \\
1048576,10 & 1.70e-05 & 2.20 & 2.87e+02 & 5.86e+04 & 5.03e+00 & 1.17e+04 \\
\bottomrule
\end{tabular}
\caption{Numerical results for the one-dimensional FIO given
  in \eqref{eqn:1D-FIO}. }
\label{tab:1D-FIO}
\end{table}

Table~\ref{tab:1D-FIO} presents two groups of numerical results.
The first group adopts $7$ Chebyshev points and the relative error
is around 8.00e-03, whereas the second group adopts $10$ Chebyshev points
and the relative error is around 1.00e-05.
Theoretically, the relative error should be independent of problem size $N$
\cite{FIO09,FIO14}.
In practice, even though the relative error increases slowly
as the size of the problem increases due to the accumulation
of the numerical error, the error stays in the same order.
The third column of the table indicates that the compression ratio
is around 2.3 for the first group and 2 for the second group.
This implies that the sweeping compression procedure in the nearly optimal
IBF compresses the factorization by a factor greater than 2,
which results in saving in both memory and application time.
The saving is greater in higher dimensional problems as we can see in
previous analysis and in the example later.
On the time scaling side, both the factorization time and the application
time strongly support the complexity analysis.
Every time we quadripule the problem size,
the factorization time increases on average by a factor of 5,
and the increasing factor decreases monotonically down to 4.
The increasing factor for the application time is on average
lower but close to 4.
The speedup factor over the direct method may catch the eye ball
of the users who are interested in the application of the FIO.

\subsection{General Fourier Integral Operator in 2D with MIBF}
\label{sub:General Fourier Integral Operator in 2D with MBF}

This section presents a numerical example to demonstrate the
efficiency of the multiscale interpolative butterfly factorization (MIBF).

We revisit a similar example in \cite{MBF},
\begin{equation}\label{eq:2D-FIO-MBF}
u(x) = \sum_{\xi\in\Omega} e^{2\pi \imath \Phi(x,\xi)}g(\xi),
\quad x\in X,
\end{equation}
with a kernel $\Phi(x,\xi)$ given by
\begin{equation}
\begin{split}
\Phi(x,\xi) =& x\cdot \xi+\sqrt{c_1^2(x)\xi_1^2+c_2^2(x)\xi_2^2},\\
c_1(x) = & (2+\sin(2\pi x_1)\sin(2\pi x_2))/32,\\
c_2(x) = & (2+\cos(2\pi x_1)\cos(2\pi x_2))/32,
\end{split}
\end{equation}
where $X$ and $\Omega$ are defined as,
\begin{equation}
  \label{eq:X}
  X = \left\{ x = \left( \frac{n_1}{N^{1/2}}, \frac{n_2}{N^{1/2}}\right),
      0 \leq  n_1,n_2 < N^{1/2}\text{ with }
  n_1, n_2 \in \mathbb{Z} \right\}
\end{equation}
and
\begin{equation}
  \label{eq:Omega}
  \Omega = \left\{ \xi = (n_1, n_2),
      - \frac{N^{1/2}}{2} \leq n_1,n_2 < \frac{N^{1/2}}{2}\text{ with }
  n_1, n_2 \in \mathbb{Z} \right\}.
\end{equation}
In the multiscale decomposition of $\Omega$, we recursively divide $\Omega$
until the center part is of size 16 by 16.

\begin{table}[htp]
\centering
\begin{tabular}{rcccccc}
\toprule
   $N,r_\epsilon$ & $\epsilon^m$ & $R_{comp}$ & $T_{Factor}(min)$
                             & $T_d(sec)$ & $T_{app}(sec)$ & $T_d/T_{app}$\\
\toprule
 $32^2$, 6 & 2.52e-03 & 2.45 & 7.63e-02 & 2.45e-01 & 7.52e-03 & 3.25e+01 \\
 $64^2$, 6 & 4.13e-03 & 2.47 & 4.00e-01 & 2.17e+00 & 3.71e-02 & 5.86e+01 \\
$128^2$, 6 & 3.11e-03 & 2.41 & 2.19e+00 & 2.11e+01 & 2.65e-01 & 7.98e+01 \\
$256^2$, 6 & 1.71e-02 & 3.08 & 1.91e+01 & 2.61e+02 & 1.19e+00 & 2.20e+02 \\
$512^2$, 6 & 5.32e-02 & 3.35 & 9.58e+01 & 4.88e+03 & 5.59e+00 & 8.74e+02 \\
\toprule
 $32^2$, 9 & 5.58e-06 & 3.79 & 1.77e-01 & 2.44e-01 & 1.17e-02 & 2.08e+01 \\
 $64^2$, 9 & 7.21e-06 & 2.96 & 1.07e+00 & 2.27e+00 & 7.68e-02 & 2.95e+01 \\
$128^2$, 9 & 6.98e-06 & 2.66 & 5.55e+00 & 2.09e+01 & 6.12e-01 & 3.41e+01 \\
$256^2$, 9 & 8.37e-06 & 3.16 & 6.34e+01 & 2.85e+02 & 8.48e+00 & 3.36e+01 \\
$512^2$, 9 & 1.23e-05 & 2.95 & 3.11e+02 & 4.79e+03 & 5.25e+01 & 9.13e+01 \\
\bottomrule
\end{tabular}
\caption{Numerical results for the two-dimensional FIO given
  in \eqref{eq:2D-FIO-MBF} by the MIBF. }
\label{tab:2D-FIO-MBF}
\end{table}

Table \ref{tab:2D-FIO-MBF} summarizes the results
of this example by the {MIBF}.
The results agree with the $\O(N\log N)$ complexity analysis.
As we double the problem size
$N$, the factorization time increases by a factor 5 on average.  Similarly,
the actual application time matches the
theoretical complexity as well.
The relative error is essentially independent of the problem size $N$
and the speedup factor is attractive.
Comparing Table~\ref{tab:1D-FIO} and Table~\ref{tab:2D-FIO-MBF},
we notice that $R_{comp}$ in two dimensions is larger than that in one-dimension.
This matches our expectation because the numerical rank by the
Chebyshev interpolation is $r_\eps^d$, which increases with the dimension
$d$. Therefore, the sweeping compression benefits more in multidimensional
problems.

\section{Conclusion and discussion}
\label{sec:conclusion}
This paper introduces an interpolative butterfly factorization as a data-sparse
approximation of complementary low-rank matrices
when their kernel functions satisfy certain analytic properties.
More precisely, it represents such an $N \times N $ dense matrix
as a product of $\O(\log N)$ sparse matrices
with nearly optimal number of entries.
The construction and application of
the interpolative butterfly factorization is highly efficient
with $O(N\log N)$ operation and memory complexity.
The prefactor of the complexity is nearly optimal in the butterfly scheme.

Since applying the sparse factors is essentially a sequence of sparse
matrix-vector multiplications with structured sparsity,
this algorithm is especially of interest in distributed parallel computing.
Based on the data distribution patten given in~\cite{FIO13},
the problem can be easily distributed in a $d$-dimensional way,
which is of great interests for extreme-scale computing.
In another word, for a problem of size $N=n^d$,
we could distribute the problem among $P=\O(N)$ processes
and achieve communication complexity,
$\O( \alpha \log p + \beta \frac{N}{P}r_0 \log P )$,
where $\alpha$ is the message latency
and $\beta$ is the per-process inverse bandwidth.
It is a promising general framework for scalable implementation
of a wide range of transforms in harmonic analysis.

{\bf Acknowledgments.}  Y. Li was partially
supported by the National Science Foundation under award DMS-1328230
and the U.S. Department of Energy's Advanced Scientific Computing
Research program under award DE-FC02-13ER26134/DE-SC0009409. 
H. Y. is partially supported by the National Science Foundation 
under grants ACI-1450280 and thank the support of the 
AMS-Simons travel award.

\bibliographystyle{abbrv}
\bibliography{ref}

\begin{thebibliography}{10}

\bibitem{PolarFFT}
A.~Averbuch, R.~Coifman, D.~Donoho, M.~Elad, and M.~Israeli.
\newblock Fast and accurate polar {F}ourier transform.
\newblock {\em Applied and Computational Harmonic Analysis}, 21(2):145 -- 167,
  2006.

\bibitem{PFFT2}
O.~Ayala and L.-P. Wang.
\newblock Parallel implementation and scalability analysis of 3{D} fast
  {F}ourier transform using 2{D} domain decomposition.
\newblock {\em Parallel Computing}, 39(1):58 -- 77, 2013.

\bibitem{symbol1}
G.~Bao and W.~Symes.
\newblock Computation of pseudo-differential operators.
\newblock {\em SIAM Journal on Scientific Computing}, 17(2):416--429, 1996.

\bibitem{bernstein2007tangent}
D.~J. Bernstein.
\newblock The tangent {FFT}.
\newblock In {\em Applied algebra, algebraic algorithms and error-correcting
  codes}, pages 291--300. Springer, 2007.

\bibitem{FIO07}
E.~Cand{\`e}s, L.~Demanet, and L.~Ying.
\newblock Fast computation of {F}ourier integral operators.
\newblock {\em SIAM J. Sci. Comput.}, 29(6):2464--2493, 2007.

\bibitem{PR1}
E.~Cand{\`e}s, X.~Li, and M.~Soltanolkotabi.
\newblock Phase retrieval via {W}irtinger flow: Theory and algorithms.
\newblock {\em Information Theory, IEEE Transactions on}, 61(4):1985--2007,
  April 2015.

\bibitem{FIO09}
E.~J. Cand{\`e}s, L.~Demanet, and L.~Ying.
\newblock A fast butterfly algorithm for the computation of {{F}ourier}
  integral operators.
\newblock {\em Multiscale Modeling and Simulation}, 7(4):1727--1750, 2009.

\bibitem{Imaging1}
L.~Demanet, M.~Ferrara, N.~Maxwell, J.~Poulson, and L.~Ying.
\newblock A butterfly algorithm for synthetic aperture radar imaging.
\newblock {\em SIAM Journal on Imaging Sciences}, 5(1):203--243, 2012.

\bibitem{symbol2}
L.~Demanet and L.~Ying.
\newblock Discrete symbol calculus.
\newblock {\em SIAM Review}, 53(1):71--104, 2011.

\bibitem{FIO12}
L.~Demanet and L.~Ying.
\newblock Fast wave computation via {{F}ourier} integral operators.
\newblock {\em Mathematics of Computation}, 81:1455--1486, 2012.

\bibitem{Yingwave}
L.~Demanet and L.~Ying.
\newblock Fast wave computation via {{F}ourier} integral operators.
\newblock {\em Math. Comput.}, 81(279), 2012.

\bibitem{Butterfly3}
B.~Engquist and L.~Ying.
\newblock Fast directional multilevel algorithms for oscillatory kernels.
\newblock {\em SIAM Journal on Scientific Computing}, 29(4):1710--1737, 2007.

\bibitem{Butterfly4}
B.~Engquist and L.~Ying.
\newblock A fast directional algorithm for high frequency acoustic scattering
  in two dimensions.
\newblock {\em Communications in Mathematical Sciences}, 7(2):327--345, 06
  2009.

\bibitem{FFT}
W.~M. Gentleman and G.~Sande.
\newblock Fast {F}ourier transforms: For fun and profit.
\newblock In {\em Proceedings of the November 7-10, 1966, Fall Joint Computer
  Conference}, AFIPS '66 (Fall), pages 563--578, New York, NY, USA, 1966. ACM.

\bibitem{NUFFT}
L.~Greengard and J.-Y. Lee.
\newblock Accelerating the {N}onuniform {F}ast {F}ourier {T}ransform.
\newblock {\em SIAM Review}, 46(3):443--454, 2004.

\bibitem{Rec1}
N.~Halko, P.~Martinsson, and J.~Tropp.
\newblock Finding structure with randomness: Probabilistic algorithms for
  constructing approximate matrix decompositions.
\newblock {\em SIAM Review}, 53(2):217--288, 2011.

\bibitem{Hu}
J.~Hu, S.~Fomel, L.~Demanet, and L.~Ying.
\newblock {A fast butterfly algorithm for generalized {Radon} transforms}.
\newblock {\em Geophysics}, 78(4):U41--U51, June 2013.

\bibitem{johnson2007modified}
S.~G. Johnson and M.~Frigo.
\newblock A modified split-radix {FFT} with fewer arithmetic operations.
\newblock {\em Signal Processing, IEEE Transactions on}, 55(1):111--119, 2007.

\bibitem{Knockaert}
L.~Knockaert.
\newblock Fast {H}ankel transform by fast sine and cosine transforms: the
  {M}ellin connection.
\newblock {\em Signal Processing, IEEE Transactions on}, 48(6):1695--1701, Jun
  2000.

\bibitem{BF}
Y.~Li, H.~Yang, E.~R. Martin, K.~L. Ho, and L.~Ying.
\newblock Butterfly {F}actorization.
\newblock {\em Multiscale Modeling \& Simulation}, 13(2):714--732, 2015.

\bibitem{MBF}
Y.~Li, H.~Yang, and L.~Ying.
\newblock Multidimensional {B}utterfly {F}actorization.
\newblock {\em arXiv:1509.07925 [math.NA]}, 2015.

\bibitem{FIO14}
Y.~Li, H.~Yang, and L.~Ying.
\newblock A multiscale butterfly aglorithm for {{F}ourier} integral operators.
\newblock {\em Multiscale Modeling and Simulation}, 13(2):614--631, 2015.

\bibitem{SP2}
W.~Liao and A.~Fannjiang.
\newblock {MUSIC} for single-snapshot spectral estimation: Stability and
  super-resolution.
\newblock {\em Applied and Computational Harmonic Analysis}, 40(1):33 -- 67,
  2016.

\bibitem{Butterfly1}
E.~Michielssen and A.~Boag.
\newblock A multilevel matrix decomposition algorithm for analyzing scattering
  from large structures.
\newblock {\em Antennas and Propagation, IEEE Transactions on},
  44(8):1086--1093, Aug 1996.

\bibitem{Butterfly2}
M.~O'Neil, F.~Woolfe, and V.~Rokhlin.
\newblock An algorithm for the rapid evaluation of special function transforms.
\newblock {\em Appl. Comput. Harmon. Anal.}, 28(2):203--226, 2010.

\bibitem{PFFT1}
D.~Pekurovsky.
\newblock P3{D}{FFT}: A framework for parallel computations of {F}ourier
  transforms in three dimensions.
\newblock {\em SIAM Journal on Scientific Computing}, 34(4):C192--C209, 2012.

\bibitem{PFFT}
D.~Pekurovsky.
\newblock {P3{D}FFT}: {A} {F}ramework for {P}arallel {C}omputations of
  {F}ourier {T}ransforms in {T}hree {D}imensions.
\newblock {\em SIAM Journal on Scientific Computing}, 34(4):C192--C209, 2012.

\bibitem{FIO13}
J.~Poulson, L.~Demanet, N.~Maxwell, and L.~Ying.
\newblock A parallel butterfly algorithm.
\newblock {\em SIAM J. Sci. Comput.}, 36(1):C49--C65, 2014.

\bibitem{SHT2}
V.~Rokhlin and M.~Tygert.
\newblock Fast algorithms for spherical harmonic expansions.
\newblock {\em SIAM J. Sci. Comput.}, 27(6):1903--1928, Dec. 2005.

\bibitem{wavemoth}
D.~S. Seljebotn.
\newblock Wavemoth-fast spherical harmonic transforms by butterfly matrix
  compression.
\newblock {\em The Astrophysical Journal Supplement Series}, 199(1):5, 2012.

\bibitem{Alex}
A.~Townsend.
\newblock A fast analysis-based discrete {H}ankel transform using asymptotic
  expansions.
\newblock {\em SIAM Journal on Numerical Analysis}, 53(4):1897--1917, 2015.

\bibitem{invRadon}
D.~O. Trad, T.~J. Ulrych, and M.~D. Sacchi.
\newblock Accurate interpolation with high-resolution time-variant {Radon}
  transforms.
\newblock {\em Geophysics}, 67(2):644--656, 2002.

\bibitem{SHT}
M.~Tygert.
\newblock Fast algorithms for spherical harmonic expansions, \{III\}.
\newblock {\em Journal of Computational Physics}, 229(18):6181 -- 6192, 2010.

\bibitem{symbol3}
H.~Yang and L.~Ying.
\newblock A fast algorithm for multilinear operators.
\newblock {\em Applied and Computational Harmonic Analysis}, 33(1):148 -- 158,
  2012.

\bibitem{invFIO}
B.~Yazici, L.~Wang, and K.~Duman.
\newblock Synthetic aperture inversion with sparsity constraints.
\newblock In {\em Electromagnetics in Advanced Applications (ICEAA), 2011
  International Conference on}, pages 1404--1407, Sept 2011.

\bibitem{Butterfly5}
L.~Ying.
\newblock Sparse {{F}ourier} transform via butterfly algorithm.
\newblock {\em SIAM J. Sci. Comput.}, 31(3):1678--1694, Feb. 2009.

\bibitem{ying20053d}
L.~Ying, L.~Demanet, and E.~Cand{\`e}s.
\newblock 3{D} discrete curvelet transform.
\newblock In {\em Optics \& Photonics 2005}, pages 591413--591413.
  International Society for Optics and Photonics, 2005.

\bibitem{Zhizhen:2014}
Z.~Zhao, Y.~Shkolnisky, and A.~Singer.
\newblock {F}ast {S}teerable {P}rincipal {C}omponent {A}nalysis.
\newblock {\em Computational Imaging, IEEE Transactions on}, to appear.

\end{thebibliography}

\end{document}